\documentclass[a4paper,11pt]{article}
\usepackage[T1]{fontenc}
\usepackage{lmodern}
\usepackage[margin=28mm,headheight=14pt]{geometry}
\usepackage{amsmath,amssymb,amsthm,mathtools}
\usepackage{comment}
\usepackage{enumitem,microtype,fancyhdr,xcolor,xurl,needspace}
\usepackage[colorlinks=true,linkcolor=blue!40!black,citecolor=blue!40!black,urlcolor=blue!40!black]{hyperref}
\hypersetup{pdftitle={Delaunay classification of positive singular solutions to the critical Lane-Emden system},pdfauthor={Xia Huang, Yahui Jiang, Feng Zhou, Xianmei Zhou}}
\allowdisplaybreaks[2]
\setlist[enumerate]{label=\textup{(\roman*)},leftmargin=2.3em,itemsep=0.25em}
\newtheorem{theorem}{Theorem}[section]
\newtheorem{proposition}[theorem]{Proposition}
\newtheorem{lemma}[theorem]{Lemma}
\newtheorem{corollary}[theorem]{Corollary}
\theoremstyle{remark}

\numberwithin{equation}{section}
\newcommand{\R}{\mathbb R}
\newcommand{\dd}{\,\mathrm d}
\newcommand{\ee}{\mathrm e}
\newcommand{\eps}{\varepsilon}
\newcommand{\one}{\mathbf1}
\newcommand{\cU}{\mathcal U}
\newcommand{\cV}{\mathcal V}
\newcommand{\cH}{\mathcal H}
\newcommand{\K}{\mathcal K}
\newcommand{\Sph}{\mathbb S^{n-1}}
\newcommand{\om}{\omega_{n-1}}

\title{\bfseries Delaunay classification of positive singular solutions\\
to the critical Lane--Emden system}
\author{Xia Huang\textsuperscript{1}\quad Yahui Jiang\textsuperscript{1}\quad
Feng Zhou\textsuperscript{1}\quad Xianmei Zhou\textsuperscript{2}\\[0.8em]
\begin{minipage}{0.92\textwidth}\centering\small
\textsuperscript{1}School of Mathematical Sciences, Key Laboratory of MEA
(Ministry of Education) \& Shanghai Key Laboratory of PMMP,\\
East China Normal University\\
Xia Huang: \texttt{xhuang@cpde.ecnu.edu.cn};
Yahui Jiang: \texttt{52275500028@stu.ecnu.edu.cn}\\
Feng Zhou: \texttt{fzhou@math.ecnu.edu.cn}\\[0.5em]
\textsuperscript{2}School of Mathematical Sciences, Zhejiang Normal University\\
Xianmei Zhou: \texttt{xmzhou@zjnu.edu.cn}
\end{minipage}}
\date{}

\begin{document}
\maketitle
\begin{abstract}
We classify positive classical solutions of the Lane--Emden system
\[
 \begin{cases}
 -\Delta u=v^p & \text{in }\R^n\setminus\{0\},\\
 -\Delta v=u^q & \text{in }\R^n\setminus\{0\},
 \end{cases}
\]
with a nonremovable singularity at $0$
(i.e., \mbox{$\lim_{x\to0}u(x)=\lim_{x\to0}v(x)=+\infty$}),
where $n\ge3$ and $p,q>0$ satisfy
\[
 \frac1{p+1}+\frac1{q+1}=\frac{n-2}{n}.
\]
Without additional growth or integrability assumptions, we establish
radial symmetry about the origin, strict radial decrease, and
global two-sided power bounds for both components.
Furthermore, every such solution can be written as
\[
 u(x)=|x|^{-\alpha}U(\log|x|),\qquad
 v(x)=|x|^{-\beta}V(\log|x|),
\]
where $\alpha=n/(q+1)$, $\beta=n/(p+1)$, and $(U,V)$ is either
the positive equilibrium $(A,B)=((\alpha\beta)^{\alpha/2},
(\alpha\beta)^{\beta/2})$ or a positive nonconstant periodic solution
of the associated Emden--Fowler system. These periodic solutions form a Delaunay
family: for every $\eps\in(0,A)$, there is a unique profile pair,
up to translation in $t=\log|x|$, with $\min_{\R}U=\eps$.
The Hamiltonian energy
provides an alternative parameter, decreasing continuously and
strictly from zero to the equilibrium energy as $\eps$ increases
from $0$ to $A$. The common least period of the two components
tends to infinity as $\eps\downarrow0$ and to the period of the
linearized oscillation at $(A,B)$ as $\eps\uparrow A$. The classification includes the case of a sublinear exponent.
\end{abstract}

\section{Introduction}

The classification of positive solutions with isolated singularities
is a central problem for elliptic equations with critical growth.
On the full punctured space, it combines radial symmetry with the
determination of all positive radial profiles. The critical scalar
equations provide the basic models.
For the second-order equation
\[
 -\Delta w=w^{\frac{n+2}{n-2}},\qquad
 w>0\quad\hbox{in }\R^n\setminus\{0\},\qquad n\ge3,
\]
Caffarelli, Gidas and Spruck~\cite{CGS} proved that
a nonremovable singularity at zero forces radial symmetry about
zero, without a prescribed decay rate at infinity. The change of
variables $w(r)=r^{-(n-2)/2}\phi(\log r)$ gives
\[
 \phi''-\frac{(n-2)^2}{4}\phi
       +\phi^{\frac{n+2}{n-2}}=0.
\]
Its phase portrait yields the singular classification: $\phi$ is
either the positive constant equilibrium or a positive nonconstant
periodic profile. The periodic profiles form the Fowler, or
Delaunay, family, uniquely determined by their minimum up to
translation; see~\cite{KMPS}. The positive homoclinic
orbit gives a smooth entire bubble. All positive smooth entire
solutions are classified in~\cite{CGS,CL}.
For the local problem in a punctured ball, Caffarelli, Gidas and
Spruck~\cite{CGS} proved asymptotic radial symmetry, and
Korevaar, Mazzeo, Pacard and Schoen~\cite{KMPS} obtained refined
asymptotics in terms of deformed Fowler profiles. These local
asymptotic results are distinct from classification on the full
punctured space.

The fourth-order counterpart is
\[
 \Delta^2 w=w^{\frac{n+4}{n-4}},\qquad
 w>0\quad\hbox{in }\R^n\setminus\{0\},\qquad n\ge5.
\]
Lin~\cite{Lin} classified the smooth positive entire solutions;
his moving-plane argument also gives radial symmetry
about zero for solutions with a nonremovable origin. The substitution
$w(r)=r^{-(n-4)/2}\psi(\log r)$ leads to
\[
 \psi^{(4)}-\frac{n^2-4n+8}{2}\psi''
       +\frac{n^2(n-4)^2}{16}\psi
       -\psi^{\frac{n+4}{n-4}}=0.
\]
This equation has a four-dimensional phase space, so the radial
classification requires additional structure. Guo, Huang, Wang
and Wei~\cite{GHWW} constructed positive periodic profiles by
variational methods for sufficiently large prescribed periods and
studied oscillation about the homogeneous singular solution.
Frank and K\"onig~\cite{FK} obtained the complete
classification: every positive solution with a nonremovable origin
has a constant or periodic profile, and the nonconstant periodic
profile is unique up to translation for each prescribed minimum.
They also proved strict radial decrease. For nonremovable compact
singular sets contained in a hyperplane, Huang, Li and
Zhou~\cite{HLZ} established reflection symmetry and monotonicity in
the normal direction under a suitable upper Minkowski dimension bound.

We consider the corresponding singular problem for the critical
Lane--Emden system
\begin{equation}\label{system}
 -\Delta u=v^p,\qquad -\Delta v=u^q,
 \qquad u,v>0\quad\hbox{in }\R^n\setminus\{0\},
\end{equation}
where
\begin{equation}\label{critical}
 n\ge3,\qquad p,q>0,\qquad
 \frac1{p+1}+\frac1{q+1}=\frac{n-2}{n}.
\end{equation}
The diagonal class $u=v$ with $p=q=(n+2)/(n-2)$ recovers the
second-order model. If $p=1$, then $n\ge5$ and
$q=(n+4)/(n-4)$; eliminating $v=-\Delta u>0$ gives the critical
biharmonic equation for $u$. The analogous reduction holds when
$q=1$.

Two difficulties arise for the general system. Radial symmetry must
be proved without assuming global integrability. Moreover, when
$p\ne q$, the different scaling weights produce a four-dimensional
Emden--Fowler system with first-order terms of opposite signs.
The phase-plane argument for the second-order scalar equation
therefore requires additional structure of the positive trajectories
before it can be applied to this coupled system.

For the associated integral system on $\R^n$, Chen, Li and
Ou~\cite{CLO} proved radial symmetry in the finite-energy class
$u\in L^{q+1}(\R^n)$, $v\in L^{p+1}(\R^n)$, when $p,q\ge1$.
Hang~\cite{Hang} established smoothness and radial symmetry about
a common center in this class throughout the admissible range,
including sublinear exponents.
Hulshof and van der Vorst~\cite{HV} determined the asymptotic
behavior of radial ground states. Frank, Kim and
Pistoia~\cite{FKP} proved their nondegeneracy modulo translations
and dilations, in the natural energy class and for perturbations
vanishing at infinity. For radial solutions of weighted critical
H\'enon--Lane--Emden systems, Musina and Sreenadh~\cite{MS}
developed variational and Emden--Fowler methods.

\Needspace{6\baselineskip}
Singular critical systems with a common Sobolev homogeneity have
also been studied. Ghergu, Kim and Shahgholian~\cite{GKS}
proved radial symmetry of nonnegative classical solutions of
$-\Delta\mathbf w=|\mathbf w|^{4/(n-2)}\mathbf w$ on
$\R^n\setminus\{0\}$ with a nonremovable origin, and analyzed
local singular behavior. The common vector nonlinearity in that
system differs from the mixed powers in~\eqref{system}, so its
symmetry theory does not cover general points of~\eqref{critical}.

Guo, Hu, Peng and Yuan~\cite{GHPY} recently classified positive
smooth entire solutions of the critical Lane--Emden system.
Here we study positive solutions with a nonremovable singularity
at the origin and establish radial symmetry and a complete
classification into homogeneous and Delaunay solutions, without
additional growth or integrability assumptions. The classification
covers the full critical hyperbola, including sublinear exponents.

Throughout, solutions of~\eqref{system}--\eqref{critical} belong to
$C^2(\R^n\setminus\{0\})^2$. The origin is \emph{nonremovable}
if the pair has no positive $C^2$ extension satisfying the system
across zero. In this setting, nonremovability is equivalent to
$\lim_{x\to0}u(x)=\lim_{x\to0}v(x)=+\infty$, by
Lemma~\ref{Newton}\textup{(iii)} and the radial monotonicity in
Theorem~\ref{symmetry}. Set
\begin{equation*}
 \alpha=\frac n{q+1},\qquad \beta=\frac n{p+1},\qquad k=n-2.
\end{equation*}
Then
\begin{equation}\label{relations}
 \alpha+\beta=k,\qquad p\beta=\alpha+2,
 \qquad q\alpha=\beta+2,\qquad pq>1,
\end{equation}
and $p,q>2/k$. The common scaling
$u_\lambda(x)=\lambda^\alpha u(\lambda x)$,
$v_\lambda(x)=\lambda^\beta v(\lambda x)$ preserves the system
and its critical integrals. For homogeneous singular solutions,
the critical densities have order $r^{-n}$: their integrals diverge
at both zero and infinity, while their integrals over
$R<|x|<2R$ are independent of $R$. A global finite-energy
hypothesis would therefore exclude these basic singular solutions.

The first theorem establishes radial symmetry, strict radial decrease,
and the critical pointwise bounds.

\begin{theorem}\label{symmetry}
Let $(u,v)$ solve~\eqref{system}--\eqref{critical} and suppose that
the origin is nonremovable. Then both components are radial about
zero and strictly decreasing in $r=|x|$. Moreover, there is a
constant $C>0$ such that
\begin{equation}\label{pointbound}
 u(x)\le C|x|^{-\alpha},\qquad v(x)\le C|x|^{-\beta}
 \qquad(x\ne0).
\end{equation}
\end{theorem}

Since $\alpha(q+1)=\beta(p+1)=n$, estimate~\eqref{pointbound}
implies the uniform annular bound
$\sup_{R>0}\int_{R<|x|<2R}(u^{q+1}+v^{p+1})\dd x<\infty$.
For radial nonincreasing components, the converse follows by
integrating over $r/2<|x|<r$.

For a radial solution, set
\begin{equation}\label{EF}
 u(r)=r^{-\alpha}U(t),\qquad v(r)=r^{-\beta}V(t),
 \qquad t=\log r.
\end{equation}
With $c=\alpha\beta$ and $d=\beta-\alpha$, the profiles satisfy
\begin{equation}\label{ODE}
 U''+dU'-cU+V^p=0,\qquad V''-dV'-cV+U^q=0.
\end{equation}
The positive equilibrium is $(A,B)$, and the conserved Hamiltonian is
\begin{align}
 A=c^{\alpha/2},\qquad &B=c^{\beta/2},\qquad cA=B^p,
 \quad cB=A^q,\notag\\
 \cH(U,U',V,V')&=U'V'-cUV+
       \frac{U^{q+1}}{q+1}+\frac{V^{p+1}}{p+1}.
       \label{Hamiltonian}
\end{align}
Let $E_*$ be the equilibrium energy, and let $\omega>0$ be the
linear oscillation frequency, corresponding to the purely imaginary
eigenvalues $\pm i\omega$ of the linearization at $(A,0,B,0)$.
Explicitly,
\begin{align}
 E_*&=\cH(A,0,B,0)=-\frac2n cAB=-\frac2n c^{n/2},\notag\\
 \omega^2&=\frac{\sqrt{(2c+d^2)^2+4(pq-1)c^2}-(2c+d^2)}2,
 \qquad \omega>0.\label{frequency}
\end{align}

Theorem~\ref{symmetry} reduces the singular problem to the study of
positive global profiles of~\eqref{ODE}. The following theorem
classifies those corresponding to a nonremovable origin and
describes the associated periodic family.

\Needspace{20\baselineskip}
\begin{theorem}\label{classification}
Under the assumptions of Theorem~\ref{symmetry}, the profiles
in~\eqref{EF} are exactly
\[
 (U,V)\equiv(A,B)
 \quad\text{or}\quad
 (U,V)(t)=(U_\eps,V_\eps)(t+\tau).
\]
For each $\eps\in(0,A)$, there is a unique positive periodic
solution $(U_\eps,V_\eps)$ of~\eqref{ODE} normalized by
$U_\eps(0)=\min_{\R}U_\eps=\eps$.
Both components have the same least period $T_\eps$, and each has
exactly one nondegenerate minimum and maximum per least period. In the periodic
case, $\eps$ and $\tau\in[0,T_\eps)$ are unique.

The energies $E(\eps)$ form a continuous strictly decreasing
bijection $(0,A)\to(E_*,0)$, and
\begin{equation}\label{periodlimits}
 T_\eps\longrightarrow\infty\quad(\eps\downarrow0),\qquad
 T_\eps\longrightarrow\frac{2\pi}{\omega}\quad(\eps\uparrow A).
\end{equation}
\end{theorem}

As an immediate consequence, $r^\alpha u(r)$ and $r^\beta v(r)$
are bounded above and away from zero: the profiles are positive
and constant or periodic. These bounds depend on the solution.
In particular, both components tend uniformly to infinity at zero.
The family may equivalently be parameterized by $\min V\in(0,B)$.
Corollary~\ref{homoclinicorbit} identifies the unique remaining
positive global radial trajectory as the ground-state homoclinic orbit.

When $p\ne q$, a nonconstant pair of profiles cannot be even
about a common point: the odd part of the first equation
in~\eqref{ODE} would force $U'=0$, and hence both profiles would
be constant. Thus the common period does not require a common
center of evenness. Periodicity in $t=\log r$ gives discrete
scaling in the original radial variable.

For the symmetry theorem, the main issue is to obtain estimates
before imposing symmetry. Sections~\ref{sec:isolation}--\ref{sec:symmetry}
first exclude a Dirac mass at the puncture and establish exact
Newton representations and source bounds. We then adapt the
comparison estimates in~\cite{GHPY}, using ground-state
nondegeneracy~\cite{FKP}, to sequences whose punctures recede to
infinity after rescaling. The resulting isolation theorem rules
out convergence to an entire ground state while the punctures
remain nonremovable. A doubling argument~\cite{PQS} turns this
isolation into the pointwise bound needed to make the truncation
errors in the reflection comparison vanish. Nonremovability fixes
the symmetry center at zero.

For the classification theorem, Sections~\ref{sec:planar}
and~\ref{sec:Delaunay} show that the data $(U,U')$ at one time
determine a positive global trajectory and that the corresponding
planar projection has a locally Lipschitz inverse away from zero.
This permits the use of planar limit-set theory together with
local families of periodic orbits. To handle trajectories
approaching zero, Section~\ref{sec:completion} uses positive Green
operators to prove that convergence to zero as $t\to-\infty$
gives a removable origin. This argument also covers sublinear
exponents, for which the vector field is not differentiable at
zero. The entire classification then identifies the remaining
nonperiodic orbit and completes the singular classification.
Throughout, $C$ denotes a positive constant whose value may change;
dependencies beyond the fixed parameters are specified when needed.

\section{Entire limits and an escaping puncture}\label{sec:isolation}
Interchanging the two components if necessary, we assume throughout
this section that
$
 p\le q.
$
By the critical relation,
$
 q\ge \frac{n+2}{n-2}>1,
$
whereas \(p\) is allowed to be sublinear. 
Let $(\cU,\cV)$ denote a positive radial ground state of
\[
-\Delta\cU=\cV^p,\qquad
-\Delta\cV=\cU^q
\qquad\text{in }\R^n,
\]
satisfying
\[
\cU\in L^{q+1}(\R^n),\qquad
\cV\in L^{p+1}(\R^n).
\]
By \cite{HV} and Theorem~1.1 of \cite{GHPY}, such a solution
is unique up to translations and critical dilations. Moreover,
$\cU$ and $\cV$ are both radially symmetric and strictly decreasing
about a common center, and vanish at infinity.
We henceforth fix the representative centered at the origin
and normalized by $\cU(0)=1$.

Set
\[
\mathfrak M
:=\int_{\R^n}\cU^{q+1}\dd x
=\int_{\R^n}\cV^{p+1}\dd x<\infty,
\qquad
\om:=|\Sph|,
\qquad
\Gamma(x):=\frac{|x|^{-k}}{k\om},
\]
and write $I_2h:=\Gamma*h$ whenever the convolution is well defined.
Here $\mathfrak M$ denotes the common critical mass of the ground
state, which is distinct from the conserved energy $E_*$ of the
cylindrical system.

We recall three results from \cite{GHPY} that will be
used in the subsequent analysis.
First, Theorem~1.1 of \cite{GHPY} shows that every
positive entire solution is obtained from the normalized
ground state $(\cU,\cV)$ by a common translation and
critical dilation.
Second, Theorem~2.1 and Lemma~2.7 of \cite{GHPY}
establish
\begin{equation}\label{groundproperties}
\cU^q\in L^1(\R^n),
\qquad
\cV(x)\asymp(1+|x|)^{-k},
\end{equation}
together with the uniform tail estimate
\begin{equation}\label{tailratio}
\eta_R:=
\sup_{x\in\R^n}
\frac{
I_2\bigl(\cU^q\one_{\{|y|\ge R\}}\bigr)(x)
}{\cV(x)}
\longrightarrow0
\qquad\text{as }R\to\infty.
\end{equation}
Finally, Lemma~2.8 of \cite{GHPY}, combined with
the ground-state nondegeneracy result of
Frank, Kim, and Pistoia~\cite{FKP}, shows that
the admissible kernel of the linearized system
is precisely the span of the translation and
dilation modes.
Since these results apply to entire solutions,
their use in the analysis of sequences of punctured
solutions requires the isolation property
established in Proposition~\ref{isolation}.
  In the compactness and limiting
arguments below, we use the standard interior \(W^{2,s}\) and Schauder
estimates, see \cite[Chapters~6 and 8--9]{GT}.

\Needspace{12\baselineskip}
\subsection{The Newton representation at a puncture}

We first establish some basic properties of positive solutions to \eqref{system} on the punctured space. The following lemma provides whole-space Newton potential representations, local source estimates, and a removability criterion, while also excluding possible Dirac masses at the origin.

\Needspace{7\baselineskip}
\begin{lemma}\label{Newton}
Let $(u,v)\in C^2(\R^n\setminus\{0\})$ be a positive classical solution to \eqref{system} satisfying \eqref{critical}. Then the following assertions hold.
\begin{enumerate}
\item The functions \(u,v,u^q\), and \(v^p\) belong to
\(L^1_{\mathrm{loc}}(\R^n)\). Moreover, we have 
\begin{equation}\label{distribution}
 -\Delta u=v^p,\qquad-\Delta v=u^q
 \quad\hbox{in }\mathcal D'(\R^n).
\end{equation}
In particular, no Dirac mass occurs at the origin.
\item The exact representations are
\begin{equation}\label{Newtonformula}
 u=I_2(v^p),\qquad v=I_2(u^q)\quad \text{for} \ x\in\R^n\setminus\{0\}.
\end{equation}
  Furthermore, there exists a constant
\(C=C(n,p,q)>0\) such that, for every \(z\in\R^n\) and \(r>0\),
\begin{equation}\label{sourcebound}
 \int_{B_r(z)}u^q\dd x\le Cr^\alpha,\qquad
 \int_{B_r(z)}v^p\dd x\le Cr^\beta.
\end{equation}
\item If either \(u\) or \(v\) is bounded in some punctured
neighborhood of the origin, then both components extend across the
origin as positive classical solutions of the system. Equivalently,
if the origin is nonremovable, then
\[
 \sup_{B_r\setminus\{0\}}u
 =
 \sup_{B_r\setminus\{0\}}v
 =\infty
 \qquad\text{for every }r>0.
\]

\end{enumerate}
\end{lemma}

\begin{proof}
Applying the
isolated-singularity theorem of Brezis and
Lions~\cite[Theorem~1]{BL} separately to the two positive
superharmonic functions \(u\) and \(v\), we obtain
$
 u,v,u^q$ and $v^p$ belong to  $ L^1_{\mathrm{loc}}(\R^n)$, 
and
\begin{equation}\label{dis}
 -\Delta u=v^p+a_0\delta_0,\qquad
 -\Delta v=u^q+b_0\delta_0
 \qquad\text{in }\mathcal D'(\R^n), 
\end{equation}
for some constants \(a_0,b_0\ge0\). 

For each \(R>0\), let \(G_R\) denote the Dirichlet Green function of
\(B_R\). 
Using the Green representation formula in $B_R$, together with \eqref{dis}, we have  
\begin{equation}\label{localGreen}
\begin{aligned}
 u(x)
 &=h_{u,R}(x)
   +\int_{B_R}G_R(x,y)v(y)^p\,\dd y
   +a_0G_R(x,0),\\
 v(x)
 &=h_{v,R}(x)
   +\int_{B_R}G_R(x,y)u(y)^q\,\dd y
   +b_0G_R(x,0),
\end{aligned}
\qquad x\in B_R\setminus\{0\},
\end{equation}
where \(h_{u,R}\) and \(h_{v,R}\) are  the harmonic functions in
\(B_R\) satisfying
\[
 h_{u,R}=u,\qquad h_{v,R}=v
 \quad\text{on }\partial B_R.
\]
In particular,
\(h_{u,R},h_{v,R}>0\) in \(B_R\) by the maximum principle.

\medskip\noindent\emph{Step 1: We show that $a_0=b_0=0$.}
Suppose that \(a_0>0\), and fix \(\rho>0\). Since
\[
 G_\rho(x,0)=\frac{1}{k\omega_{n-1}}|x|^{-k}+O(1)
 \qquad\text{as }x\to0,
\]
together with the first identity in \eqref{localGreen}, we have 
\begin{equation}\label{atomloweru}
 u(x)\ge a_0G_\rho(x,0)\ge c|x|^{-k}
 \qquad (0<|x|<r_0)
\end{equation}
for some \(r_0>0\). If \(kq\ge n\), then
\(u^q(x)\ge c|x|^{-kq}\), which is impossible since
\(u^q\in L^1_{\mathrm{loc}}(\R^n)\). Notice that \(q>2/k\), it remains only to consider
\[
 2<kq<n.
\]

Let \(0<r<r_0/3\), take \(|x|=r\), and set 
\(A_r:=\{y:2r<|y|<3r\}\). Taking \(r_0\) sufficiently small, if necessary,
\(x,y\in B_{\rho/2}\) for \(y\in A_r\). Moreover,
\(r\le |x-y|\le4r\). Since
\(G_\rho(x,y)=\Gamma(x-y)-H_\rho(x,y)\), with \(H_\rho\) bounded on
\(B_{\rho/2}\times B_{\rho/2}\), we have
\[
 G_\rho(x,y)\ge cr^{-k},
 \qquad y\in A_r.
\]
 Using the second identity in
\eqref{localGreen} and \eqref{atomloweru}, we  yield
\begin{align*}
 v(x)
 \ge \int_{A_r}G_\rho(x,y)u(y)^q\,\dd y
 \ge cr^{-k}\int_{A_r}|y|^{-kq}\,\dd y
  \ge cr^{2-kq},
\end{align*}
and hence
\[
 v(x)^p\ge c|x|^{-p(kq-2)}
 \qquad (0<|x|<r_0).
\]
It follows from \eqref{relations} that 
$
 p(kq-2)=n+2+2q>n,
$
so \(v^p\notin L^1_{\mathrm{loc}}(\R^n)\), another contradiction.
Thus \(a_0=0\). Interchanging \((u,q)\) and \((v,p)\), we easily get 
\(b_0=0\), and \eqref{distribution} follows.

\medskip\noindent\emph{Step 2: Exact Newton representations and source bounds.}
Since  \(a_0=b_0=0\), we obtain from \eqref{localGreen} that, for every \(R>0\) and
\(x\in B_R\setminus\{0\}\),
\[
 \int_{B_R}G_R(x,y)v(y)^p\,\dd y\le u(x),\qquad
 \int_{B_R}G_R(x,y)u(y)^q\,\dd y\le v(x).
\]
For fixed \(x\ne y\), the Green function \(G_R(x,y)\) increase to 
 \(\Gamma(x-y)\) as \(R\to\infty\). Hence, for \(x\in\R^n\setminus\{0\}\), by the monotone
convergence theorem, we have
\[
 I_2(v^p)\le u,\qquad I_2(u^q)\le v.
\]
In particular, since \(u,v\in L^1_{\mathrm{loc}}(\R^n)\), it follows that \(I_2(v^p)\) and \(I_2(u^q)\)
belong to \(L^1_{\mathrm{loc}}(\R^n)\).

Define
\[
 w_u:=u-I_2(v^p),\qquad
 w_v:=v-I_2(u^q).
\]
  Then 
\[
 -\Delta w_u=0,\qquad -\Delta w_v=0
 \quad\text{in }\mathcal D'(\R^n),
\]
with 
\(w_u,w_v\ge0\) almost everywhere. By Weyl's lemma and Liouville's theorem for harmonic functions bounded from below, there exist constants \(A,B\ge0\) such that
\[
 u=I_2(v^p)+A,\qquad
 v=I_2(u^q)+B
 \quad\text{a.e. in }\R^n.
\]
We claim that \(A=B=0\). Indeed, if \(A>0\), then
\(u^q\ge A^q\), and therefore
\[
 I_2(u^q)(x)
 \ge A^q\int_{\R^n}\Gamma(x-y)\,\dd y
 =\infty.
\]
This contradicts
\[
 v(x)=I_2(u^q)(x)+B<\infty, \quad x\in \R^n\setminus\{0\}.
\]
Thus \(A=0\). Similarly, we have 
\(B=0\), and hence
\[
 u=I_2(v^p),\qquad v=I_2(u^q).
\]

It remains to prove \eqref{sourcebound}. Fix \(z\in\R^n\) and
\(r>0\), and define
\[
 F:=\int_{B_r(z)}u^q\,\dd x,
 \qquad
 G:=\int_{B_r(z)}v^p\,\dd x.
\]
If \(x,y\in B_r(z)\), then \(|x-y|\le2r\), and hence
\[
 \Gamma(x-y)
 =\frac{1}{k\omega_{n-1}}|x-y|^{-k}
 \ge cr^{-k}.
\]
By the Newton representations,  we have 
\[
 u(x)\ge cr^{-k}G,
 \qquad
 v(x)\ge cr^{-k}F. \quad\text{for a.e. }x\in B_r(z).
\]
Raising the first inequality to the power \(q\), the second to the
power \(p\), and integrating over \(B_r(z)\), we obtain
\begin{equation}\label{coupledmass}
 F\ge cr^{n-kq}G^q,
 \qquad
 G\ge cr^{n-kp}F^p.
\end{equation}
Substituting the second inequality in Eq.~\eqref{coupledmass} into
the first,  we yield
\[
 F^{pq-1}
 \le Cr^{\,kq+kpq-n(q+1)}.
\]
Similarly, there holds 
\[
 G^{pq-1}
 \le Cr^{\,kp+kpq-n(p+1)}.
\]
By \eqref{critical}, we know that 
\[
 kq+kpq-n(q+1)=2(p+1),
 \qquad
 kp+kpq-n(p+1)=2(q+1).
\]
and 
\[
 \alpha=\frac{2(p+1)}{pq-1},
 \qquad
 \beta=\frac{2(q+1)}{pq-1},
\]
we conclude that
\[
 \int_{B_r(z)}u^q\,\dd x\le Cr^\alpha,
 \qquad
 \int_{B_r(z)}v^p\,\dd x\le Cr^\beta,
\]
where  \(C>0\) depending only on \(n,p\), and \(q\).

\medskip\noindent\emph{Step 3: Boundedness implies removability.}
Suppose that \(u\) is bounded in a punctured neighborhood of the
origin. Choose \(\rho>0\) sufficiently small such that
\[
 \|u\|_{L^\infty(B_{2\rho}\setminus\{0\})}<\infty.
\]
Then \(u^q\in L^\infty(B_{2\rho})\). We claim that \(v\) is also
bounded near the origin.

Fix \(x_*\in B_\rho\setminus\{0\}\). For
\(x\in B_\rho\setminus\{0\}\), we decompose
\[
 v(x)
 =
 \int_{B_{2\rho}}\Gamma(x-y)u(y)^q\,\dd y
 +
 \int_{\R^n\setminus B_{2\rho}}
       \Gamma(x-y)u(y)^q\,\dd y:=v_{1}(x)+v_{2}(x).
\]
Since
\(u^q\in L^\infty(B_{2\rho})\) and
\[
 \sup_{x\in B_\rho}
 \int_{B_{2\rho}}|x-y|^{-k}\,\dd y
 \le C\rho^2, 
\]
the local term $v_{1}$ is  uniformly bounded in \(B_\rho\setminus\{0\}\).
On the other hand, for \(x\in B_\rho\) and \(y\notin B_{2\rho}\), we have  
$
 \Gamma(x-y)\le C\Gamma(x_*-y),
$
and therefore
\[
\begin{aligned}
 v_{2}(x)
 \le
 C\int_{\R^n\setminus B_{2\rho}}
 \Gamma(x_*-y)u(y)^q\,\dd y
 \le Cv(x_*)<\infty.
\end{aligned}
\]
Thus \(v\) is bounded in \(B_\rho\setminus\{0\}\). Consequently,
both \(u^q\) and \(v^p\) are bounded near the origin.

Since the equations \eqref{distribution} hold distributionally across the origin and their right-hand sides are bounded, standard interior elliptic regularity and Sobolev embedding yield
\[
u,v\in W^{2,s}(B_{\rho/2})
\hookrightarrow C^{1,\gamma}(B_{\rho/2})
\]
for any \(s>n\), where \(\gamma=1-n/s\in(0,1)\). In particular, $u$ and $v$ extend continuously and nonnegatively across the origin.

Since $t\mapsto t^a$ is locally H\"older continuous on
$[0,\infty)$ for every $a>0$,  both $u^q$ and $v^p$ are locally H\"older continuous in $B_{\rho/2}$.
By the standard Schauder estimates, we yield
\[
u,v\in C^{2,\theta}(B_{\rho/4})
\]
for some $\theta\in(0,1)$. Consequently, $(u,v)$
satisfies
\[
-\Delta u=v^p,\qquad -\Delta v=u^q 
\]
holds pointwise in \(B_{\rho/4}\).
Thus $(u,v)$ can be extended across the origin as a classical solution.

Finally, by the strong minimum principle, we  show that
\(u(0),v(0)>0\). Indeed, suppose that  \(u(0)=0\). Since 
\[
 -\Delta u=v^p\ge0,\qquad -\Delta v=u^q\ge0 \quad \text{in} \ B_{\rho/4},
\]
\(u\) would attain its
minimum at an interior point and hence would be constant, contradicting
\(u>0\) in \(B_{\rho/4}\setminus\{0\}\). The same argument applies to
\(v\). Thus \(u(0),v(0)>0\).
Therefore the origin is removable.
\end{proof}

\subsection{Isolation of the entire limit}

The following proposition extends the ground-state isolation result
of~\cite[Theorem~2.5]{GHPY} to sequences of solutions on punctured
spaces, with the punctures tending to infinity.

\begin{proposition}\label{isolation}
Assume \eqref{critical} holds. Let $a_j\in\R^n$ satisfy
	$|a_j|\to\infty$, and let $(u_j,v_j)$ be positive classical
	solutions of
	\begin{equation}\label{aj}
    -\Delta u_j=v_j^p,\qquad -\Delta v_j=u_j^q
	\quad\text{in }\R^n\setminus\{a_j\}.
  \end{equation}
	Suppose that
	\[
	(u_j,v_j)\rightarrow(U_*,V_*)
	\quad\text{in }C^2_{\mathrm{loc}}(\R^n).
	\]
	for some positive entire solution $(U_*,V_*)$
	of the system. Then, for all sufficiently large $j$,
	$(u_j,v_j)$ extends to a positive classical solution on $\R^n$.
\end{proposition}

\begin{proof}
By \cite[Theorem~1.1]{GHPY}, the limit $(U_*,V_*)$ is a translation and critical dilation of a positive radial ground state. Applying a suitable translation and critical dilation to $(u_j,v_j)$, still denoted by $(u_j,v_j)$, we may assume that $(u_j,v_j)\to(\cU,\cV)$ in $C^2_{\mathrm{loc}}(\R^n)$, where $(\cU,\cV)$ is radially symmetric about the origin and normalized by $\cU(0)=1$. Under this change of variables, we also have $|a_j|\to\infty$.

Since $\cU$ is radial and strictly decreasing with respect to
$r=|x|$, and $-\Delta\cU=\cV^p$ in $\R^n$, we have
\begin{equation}\label{eq:radial-properties}
\nabla\cU(0)=0,\qquad
D^2\cU(0)=-\frac{\cV(0)^p}{n}\mathrm{Id}.
\end{equation}
Thus, the origin is a strict nondegenerate maximum point of
$\cU$. Choose a sufficiently small ball $B_\delta$ centered
at the origin. By local $C^2$ convergence, $u_j$ attains
an interior maximum at some $x_j\in B_\delta$ for all
sufficiently large $j$, with
\[
x_j\to0,\qquad \nabla u_j(x_j)=0,\qquad
u_j(x_j)\to\cU(0)=1.
\]
Set $\rho_j=u_j(x_j)^{-1/\alpha}$ and define
\[
\widetilde u_j(x)
=\rho_j^\alpha u_j(x_j+\rho_j x),
\qquad
\widetilde v_j(x)
=\rho_j^\beta v_j(x_j+\rho_j x).
\]
Then
\[
\widetilde u_j(0)=1,\qquad
\nabla\widetilde u_j(0)=0.
\]
Since $x_j\to0$ and $\rho_j\to1$, we have
\[
(\widetilde u_j,\widetilde v_j)\to(\cU,\cV)
\quad\text{in }C^2_{\mathrm{loc}}(\R^n).
\]
Furthermore, the corresponding puncture points $(a_j-x_j)/\rho_j\to\infty$ as $j\to\infty$.
Relabelling this sequence, we may assume that
\begin{equation}\label{gauge}
	u_j(0)=1,\qquad \nabla u_j(0)=0.
\end{equation}

By \eqref{groundproperties} and \eqref{sourcebound}, we see
\[
 S_j:=\int_{\R^{n}}\cV(x)v_j^p(x)\,dx
 \le C\left(1+\sum_{m\ge0}2^{-m(k-\beta)}\right)<\infty.
\]
It follows from Lemma~\ref{Newton} and Tonelli's theorem that
\[
S_j=\int_{\R^n}\cV(x)v_j^p(x)\,dx=\int_{\R^n}\int_{\R^n}\Gamma(x-y)\cU^q(y)v_j^p(x)\,dydx=\int_{\R^n}\cU^q(x)u_j(x)\,dx.
\]
Moreover, by H\"older's inequality and \eqref{Newtonformula}, we have 
\[
 I_2(\cU^{q-1}u_j)
 \le v_j^{1/q}\cV^{(q-1)/q},
\]
and hence
\begin{equation}\label{eq:S_j1}
S_j=\int_{\R^n}\cV^{p+1}\biggl(\frac{v_j}{\cV}\biggr)^{p}\ge\int_{\R^n}\cV^{p+1}
\left(\frac{I_2(\cU^{q-1}u_j)}{\cV}\right)^{pq}.
\end{equation}
Since $pq>1$, H\"older's inequality yields
\begin{equation}\label{eq:S_j2}
\left(\int_{\R^n}\cV^pI_2(\cU^{q-1}u_j)\right)^{pq}\le\left(\int_{\R^n}\cV^{p+1}\right)^{pq-1}\int_{\R^n}\cV^{p+1}\left(\frac{I_2(\cU^{q-1}u_j)}{\cV}\right)^{pq}.
\end{equation}
Combining \eqref{eq:S_j1} with \eqref{eq:S_j2}, we have
\[
\begin{aligned}
 S_j\ge\mathfrak M^{1-pq}
   \left(\int_{\R^n}\cV^pI_2(\cU^{q-1}u_j)\right)^{pq}
  =\mathfrak M^{1-pq}S_j^{pq},
 \end{aligned}
\]
which implies 
\begin{equation}\label{crossintegral}
 0<S_j=\int_{\R^n}\cU^qu_j=\int_{\R^n}\cV v_j^p\le\mathfrak M.
\end{equation}

Since $u_j\to\cU$ in $C^2_{\mathrm{loc}}(\R^n)$
and $\cU>0$, for each fixed $R>0$ and all sufficiently large $j$, we have $u_j\ge(1-o_j(1))\cU$ in $B_R$.
Together with \eqref{tailratio}, we obtain, for every $x\in\R^{n}\setminus\{a_j\}$,
\[
\begin{aligned}
	v_{j}(x)=I_2(u_j^q)(x)\ge I_2(u_j^q\mathbf1_{B_R})(x)\ge(1-o_j(1))^qI_2(\cU^q\mathbf1_{B_R})(x)\ge(1-o_j(1))^q(1-\eta_R)\cV(x).
\end{aligned}
\]
Consequently, for every $x\in\R^{n}\setminus\{a_j\}$,
\[
u_j(x)=I_2(v_j^p)(x)\ge(1-o_j(1))^{pq}(1-\eta_R)^pI_2(\cV^p)(x)=(1-o_j(1))^{pq}(1-\eta_R)^p\cU(x).
\]
Set
\begin{equation*}
	\varepsilon_j:=\max\left\{
	\sup_{x\ne a_j}\frac{(\cU-u_j)_+(x)}{\cU(x)},
	\sup_{x\ne a_j}\frac{(\cV-v_j)_+(x)}{\cV(x)}
	\right\}.
\end{equation*}
Then for each fixed $R>0$,
\[
\limsup_{j\to\infty}\varepsilon_j
\le \max\left\{
1-(1-\eta_R)^p,\eta_R
\right\}.
\]
Letting $R\to\infty$ and using \eqref{tailratio},
we conclude that $\varepsilon_j\to0$. 

We claim that $\varepsilon_j=0$ for all sufficiently
large $j$. Suppose by contradiction that there exists a
subsequence, still denoted by $\{\varepsilon_j\}$, such that
$\varepsilon_j>0$ for all $j\ge1$. 
Set
\[
\phi_j=\frac{u_j-\cU}{\varepsilon_j},
\qquad
\psi_j=\frac{v_j-\cV}{\varepsilon_j}.
\]
By the definition of $\varepsilon_j$, we have
\[
u_j\ge(1-\varepsilon_j)\cU,
\qquad
v_j\ge(1-\varepsilon_j)\cV.
\]
Then $\phi_j\ge-\cU$ and $\psi_j\ge-\cV$.
Moreover, by \eqref{crossintegral}, we have
\[
	\int_{\R^n}\cU^q(u_j-\cU)_+\,dx\le\int_{\R^n}\cU^q(\cU-u_j)_+\,dx\le\varepsilon_j\mathfrak M.
\]
Consequently,
\begin{equation}\label{linearweights1}
\int_{\R^n}\cU^q|\phi_j|\,dx
\le2\mathfrak M.
\end{equation}
Since $\varepsilon_j\to0$, for all sufficiently
large $j$ we also have
\[
(\cV^p-v_j^p)_+
\le C_p\varepsilon_j\cV^p,
\]
where $C_p=\max\{1,p\}$.
Using \eqref{crossintegral} again, we obtain
\begin{equation}\label{linearweights2}
	\int_{\R^n}\cV
	\left|\frac{v_j^p-\cV^p}{\varepsilon_j}\right|dx
	\le2C_p\mathfrak M.
\end{equation}
A direct calculation shows that 
\[
-\Delta\phi_j=A_j\psi_j,
\qquad
-\Delta\psi_j=D_j\phi_j \qquad\text{in }\R^n\setminus\{a_j\},
\]
where
\[
	A_j=p\int_0^1
	\bigl(t v_j+(1-t)\cV\bigr)^{p-1}\,dt,\quad
	D_j=q\int_0^1
	\bigl(t u_j+(1-t)\cU\bigr)^{q-1}\,dt.
\]
The local $C^2$ convergence and positivity of
$\cU,\cV$ imply that
\[
A_j\to p\cV^{p-1},
\qquad
D_j\to q\cU^{q-1}
\quad\text{in }C^2_{\mathrm{loc}}(\R^n).
\]
Consequently, \eqref{linearweights1}--\eqref{linearweights2}
yield uniform $L^1_{\mathrm{loc}}(\R^n)$ bounds
for $\phi_j$ and $\psi_j$.

Set $w_j=|\phi_j|+|\psi_j|$. Then $w_j\in L^1_{\mathrm{loc}}(\R^n)$. Since $\phi_j,\psi_j\in
C^2(\R^n\setminus\{a_j\})$, Kato's inequality \cite[Lemma~A]{Kato} and the local uniform
boundedness of $A_j,D_j$ yield that for each $R>0$,
\[
-\Delta w_j\le C_Rw_j
\quad\text{in }B_{2R}
\]
for all sufficiently large $j$.
Since
\[
\Delta_{x,t}\bigl(w_j(x)\cosh(\sqrt{C_R}\,t)\bigr)
=\cosh(\sqrt{C_R}\,t)(\Delta w_j+C_Rw_j)
\ge0,
\]
the function $w_j(x)\cosh(\sqrt{C_R}\,t)$
is nonnegative and subharmonic in
$B_{2R}\times\R$.
Thus, the mean value inequality gives, for
every $x\in B_R$,
\[
\begin{aligned}
	w_j(x)\le\frac1{|B_R^{n+1}|}
	\int_{B_R^{n+1}(x,0)}
	w_j(y)\cosh(\sqrt{C_R}\,t)\,dy\,dt\le\frac{2R\cosh(\sqrt{C_R}\,R)}
	{|B_R^{n+1}|}
	\int_{B_{2R}}w_j(y)\,dy<C_R'.
\end{aligned}
\]
Consequently, $\phi_j$ and $\psi_j$ are
uniformly bounded in $L^\infty_{\mathrm{loc}}(\R^n)$. By standard interior elliptic regularity and
Schauder estimates, up to a subsequence, $
(\phi_j,\psi_j)\to(\phi,\psi)$ in $C^2_{\mathrm{loc}}(\R^n)$, and
\[
-\Delta\phi=p\cV^{p-1}\psi,
\qquad
-\Delta\psi=q\cU^{q-1}\phi
\quad\text{in }\R^n.
\]
From the lower bounds of $\phi_j, \psi_j$, we get $\phi\ge-\cU$ and $\psi\ge-\cV$. Moreover, by Fatou's lemma, we have
\[
\int_{\R^n}\cU^q|\phi|\,dx
\le2\mathfrak M,
\qquad
p\int_{\R^n}\cV^p|\psi|\,dx
\le2C_p\mathfrak M.
\]

We now prove $(\phi,\psi)\neq(0,0)$. Set
\[
h_j=(\cU-u_j)_+,\qquad
z_j=(\cV-v_j)_+.
\]
Then
\[
\begin{aligned}
	z_j
	&=\bigl(I_2(\cU^q-u_j^q)\bigr)_+
	\le I_2\bigl((\cU^q-u_j^q)_+\bigr)
	\le qI_2(\cU^{q-1}h_j),\\
	h_j
	&=\bigl(I_2(\cV^p-v_j^p)\bigr)_+
	\le I_2\bigl((\cV^p-v_j^p)_+\bigr)
	\le C_pI_2(\cV^{p-1}z_j).
\end{aligned}
\]
Suppose that $(\phi,\psi)=(0,0)$. Then $\phi_j\to0$ in $C_{\mathrm{loc}}^2(\R^n)$, and so $h_j/(\varepsilon_j\cU)\to 0$ uniformly on every fixed ball.
Since $h_j\le\varepsilon_j\cU$, by \eqref{tailratio}, we have for each $R>0$,
\[
\begin{aligned}
\frac{z_j(x)}{\varepsilon_j\cV(x)}\le q\left\|\frac{h_j}{\varepsilon_j\cU}\right\|_{L^{\infty}(B_R)}\frac{I_2(\cU^q\mathbf1_{B_R})(x)}{\cV(x)}+q\frac{I_2(\cU^{q}\mathbf1_{\R^{n}\setminus{B_R}})(x)}{\cV(x)}
\le q\left\|
\frac{h_j}{\varepsilon_j\cU}
\right\|_{L^\infty(B_R)}
+q\eta_R.
\end{aligned}
\]
Taking the supremum over $x\ne a_j$ and then
letting $j\to\infty$, we obtain
\[
\limsup_{j\to\infty}
\left\|\frac{z_j}{\varepsilon_j\cV}\right\|_\infty
\le q\eta_R.
\]
Letting $R\to\infty$ and using
\eqref{tailratio}, we conclude that
\[
\left\|\frac{z_j}{\varepsilon_j\cV}\right\|_\infty
\rightarrow0.
\]
Similarly, we have
\[
\left\|\frac{h_j}{\varepsilon_j\cU}\right\|_{L^\infty}
\le C_p
\left\|\frac{z_j}{\varepsilon_j\cV}\right\|_{L^\infty}
\rightarrow0.
\]
This is a contradiction, since
\[
\max\left\{
\left\|\frac{h_j}{\varepsilon_j\cU}\right\|_\infty,
\left\|\frac{z_j}{\varepsilon_j\cV}\right\|_\infty
\right\}=1.
\]
Thus $(\phi,\psi)\neq(0,0)$.

On the other hand, by \cite[Lemma~2.8]{GHPY}, there exists a constant $C>0$ such that
\[
|\phi|\le C\cU,\qquad
|\psi|\le C\cV,
\]
which, together with \cite[Theorem~2]{HV}, implies $(\phi,\psi)\to(0,0)$ as $|x|\to\infty$. It then follows from \cite[Corollary~3]{FKP} that
\[
	(\phi,\psi)
	=b_0\bigl(
	\alpha\cU+x\cdot\nabla\cU,\,
	\beta\cV+x\cdot\nabla\cV
	\bigr)+\sum_{i=1}^n b_i
	(\partial_i\cU,\partial_i\cV).
\]
By \eqref{gauge} and the
local $C^2$ convergence of $\phi_j$, we obtain $\phi(0)=0$ and $\nabla\phi(0)=0$.
Combining this with \eqref{eq:radial-properties}, we obtain $b_0=b_1=\cdots=b_n=0$.
Hence $(\phi,\psi)=(0,0)$, which is a contradiction.

Therefore, $\varepsilon_j=0$ for all
sufficiently large $j$. Then $u_j\ge\cU$. Combining this with
\eqref{crossintegral}, we obtain
\[
0\le
\int_{\R^n}\cU^q(u_j-\cU)\,dx
=S_j-\mathfrak M\le0,
\]
which implies
$u_j=\cU$, and hence $v_j=\cV$. Reversing the changes of variables completes the proof.
\end{proof}

\section{Radial symmetry}\label{sec:symmetry}

In this section, we prove Theorem~\ref{symmetry}. We first establish
the critical pointwise and annular bounds for positive singular
solutions to \eqref{system}. We then prove that both $u$ and $v$ are invariant under
reflection across every hyperplane through the origin, which yields
their radial symmetry. Finally, the radial equations imply that $u$ and $v$  are strictly decreasing in the radial variable.

\subsection{Critical a priori estimates}

\begin{lemma}\label{critical-estimates}
Let \((u,v)\in C^2(\R^n\setminus\{0\})\) be a positive solution of
\eqref{system} satisfying \eqref{critical}. If the origin is
nonremovable, then there exists \(C>0\) such that
\begin{equation}\label{bound}
 u(x)\le C|x|^{-\alpha},
 \qquad
 v(x)\le C|x|^{-\beta}
 \qquad (x\ne0).
\end{equation}
Consequently,
\begin{equation}\label{intbound}
 \sup_{R>0}
 \int_{R<|x|<2R}
 \bigl(u^{q+1}+v^{p+1}\bigr)\,\dd x<\infty.
\end{equation}
\end{lemma}

\begin{proof}
Define
\[
 M(x):=u(x)^{1/\alpha}+v(x)^{1/\beta}.
\]
Suppose, to the contrary, that
\[
 \sup_{x\ne0}|x|M(x)=\infty.
\]
Then there exists a sequence \(\left\{ x_j \right\}\subset\R^n\setminus\{0\}\) 
 such that
\[
 M(x_j)|x_j|\ge2j.
\]

Applying the doubling lemma \cite[Lemma~5.1]{PQS}, there exists
\(\left\{ y_j \right\}\subset\R^n\setminus\{0\}\) such that
\[
M(y_j)\ge M(x_j),
\qquad
M(y_j)|y_j|>2j,
\]
and
\[
M(z)\le2M(y_j)
\qquad\text{whenever}\qquad
|z-y_j|\le\frac{j}{M(y_j)}.
\]
Set
$
\lambda_j:=M(y_j)^{-1},
$
and define the rescaled functions
\[
\begin{aligned}
u_j(x)&:=\lambda_j^\alpha
        u(y_j+\lambda_jx),\\
v_j(x)&:=\lambda_j^\beta
        v(y_j+\lambda_jx),
\end{aligned}
\qquad \text{for}\ x\in B_j.
\]
Since
$
j\lambda_j
=\frac{j}{M(y_j)}
<\frac{|y_j|}{2},
$
the ball \(B_{j\lambda_j}(y_j)\) does not contain the origin.
Consequently, \((u_j,v_j)\) is well-defined in \(B_j\) and satisfies
\[
-\Delta u_j=v_j^p,
\qquad
-\Delta v_j=u_j^q
\qquad\text{in }B_j.
\]
Moreover,
\[
u_j(0)^{1/\alpha}+v_j(0)^{1/\beta}
=\lambda_jM(y_j)=1,
\]
and
\[
\begin{aligned}
u_j(x)^{1/\alpha}+v_j(x)^{1/\beta}
&=\lambda_jM(y_j+\lambda_jx)
\le2 \qquad \text{for}\  x\in B_j.
\end{aligned}
\]

Fix \(R>0\). For all sufficiently large \(j\), we have
\(B_{2R}\subset B_j\), and the functions \(u_j\), \(v_j\), \(u_j^q\), and  \(v_j^p\) are uniformly
bounded in \(B_{2R}\). 
Arguing as in Step~3 of the proof of Lemma~\ref{Newton}, 
 by the interior \(W^{2,s}\) estimates and Schauder estimates, we obtain
\[
\|u_j\|_{C^{2,\eta}(B_R)}
+\|v_j\|_{C^{2,\eta}(B_R)}
\le C_R
\]
for some \(\eta\in(0,1)\), where \(C_R\) is independent of \(j\).
By a diagonal argument, up to a subsequence,
\[
(u_j,v_j)\longrightarrow(U_*,V_*)
\quad\text{in }C^2_{\mathrm{loc}}(\R^n), \qquad \text{as} 
\ j\rightarrow \infty.
\]
Therefore, we have 
\begin{equation}\label{LE}
-\Delta U_*=V_*^p,\qquad
-\Delta V_*=U_*^q
\qquad\text{in }\R^n.
\end{equation}
Moreover,
\begin{equation}\label{normlize}
  U_*(0)^{1/\alpha}+V_*(0)^{1/\beta}=1,
\qquad
U_*^{1/\alpha}+V_*^{1/\beta}\le2
\quad\text{in }\R^n.
\end{equation}
Thus \((U_*,V_*)\) is a bounded, nonnegative, nontrivial entire
solution to \eqref{LE}. Moreover, we have $U_*\not\equiv0$ and $V_*\not\equiv0$. Indeed, if \(U_*\equiv0\), then \(V_*\equiv0\), contradicting \eqref{normlize}. Then it follows from the strong
minimum principle that 
$$
U_*>0,\qquad V_*>0
\qquad\text{in }\R^n.
$$
By the classification results in~\cite[Theorem~1.1]{GHPY}, the bounded positive entire limit
\((U_*,V_*)\) is a common translate and critical dilate of the
ground state. Notice that  $(u_j,v_j)$ are  positive classical
	solutions of \eqref{aj} with $a_j=-\frac{y_{j}}{\lambda_{j}}$. Since $|a_{j}|\rightarrow \infty$, we know from Proposition~\ref{isolation} that \(a_j\) is removable for all
sufficiently large \(j\). Scaling back would make the origin
removable for \((u,v)\), contrary to the hypothesis. Hence
\[
 \sup_{x\ne0}|x|M(x)<\infty,
\]
which implies \eqref{bound}.

Finally, since \(\alpha(q+1)=\beta(p+1)=n\), for every \(R>0\), we have 
\[
\begin{aligned}
 \int_{R<|x|<2R}
 \bigl(u^{q+1}+v^{p+1}\bigr)\,\dd x
 &\le C\int_R^{2R}r^{-1}\,\dd r
 \le C,
\end{aligned}
\]
which proves \eqref{intbound}.
\end{proof}

\subsection{Centering the reflection planes}

For each \(e\in\Sph\), let
\[
\sigma_e(x):=x-2(x\cdot e)e
\]
be the reflection across the hyperplane
\[
e^\perp:=\{x\in\R^n:x\cdot e=0\}.
\]

The next lemma proves that both
$u$ and $v$ are invariant under reflection across every hyperplane
through the origin.
 Its proof adapts the
interval comparison identity of~\cite[Lemma~3.1]{GHPY}; the
truncation and centering arguments are included because the solution
is singular at the origin.

\begin{lemma}\label{centered-reflections}
Under the assumptions of Lemma~\ref{critical-estimates}, for every
\(e\in\Sph\),
\[
 u\circ\sigma_e=u,
 \qquad
 v\circ\sigma_e=v
 \quad\text{in }\R^n\setminus\{0\}.
\]
\end{lemma}

\begin{proof}
Write
\[
 f:=u^q,
 \qquad
 g:=v^p.
\]
By Lemma~\ref{critical-estimates} and the scaling identities
$q\alpha=\beta+2$ and $p\beta=\alpha+2$, we have
\begin{equation}\label{sourcepointwise}
f(x)\le C|x|^{-\beta-2},
\qquad
g(x)\le C|x|^{-\alpha-2}.
\end{equation}
Since $\alpha,\beta<n-2$, both right-hand sides are
locally integrable near the origin.

\medskip\noindent\emph{Step 1: Cutoff estimates.}
For \(R>0\), define
\[
T_R:=
\iint_{\R^n\times\R^n}
f(x)g(y)|x-y|^{1-n}
|\one_{B_R}(x)-\one_{B_R}(y)|\,\dd x\dd y.
\]

We claim that
\begin{equation}\label{cutofferror}
T_R+
\int_{\partial B_R}
\bigl(u^{q+1}+v^{p+1}\bigr)\,\dd S
\le \frac{C}{R}.
\end{equation}
Set
$
a:=\beta+2$ and 
$b:=\alpha+2$.
Since \(0<\alpha,\beta<n-2\) and
\(\alpha+\beta=n-2\), we have
$2<a,b<n$ and 
$
a+b=n+2.
$
Notice that \(|\one_{B_R}(x)-\one_{B_R}(y)|\) vanishes unless exactly one of
\(x\) and \(y\) belongs to \(B_R\).  We obtain
\[
\begin{aligned}
T_R
={}&
\int_{|x|<R}\int_{|y|>R}
\bigl(f(x)g(y)+g(x)f(y)\bigr)
|x-y|^{1-n}\,\dd y\dd x\\
\le{}&
C\int_{|x|<R}\int_{|y|>R}
\Bigl(
 |x|^{-a}|y|^{-b}
 +|x|^{-b}|y|^{-a}
\Bigr)
|x-y|^{1-n}\,\dd y\dd x.
\end{aligned}
\]

Under the change of variables \(x=RX\) and \(y=RY\), we obtain 
\[
T_R
\le
CR^{-1}
\bigl(J(a,b)+J(b,a)\bigr),
\]
where
\[
J(a,b):=
\int_{|X|<1}\int_{|Y|>1}
|X|^{-a}|Y|^{-b}|X-Y|^{1-n}\,\dd Y\dd X.
\]
It remains to verify that \(J(a,b)<\infty\).

Let
$
\mathcal D
:=
\{(X,Y)\in\R^n\times\R^n:|X|<1<|Y|\}
$
and decompose
$
\mathcal D
=
\mathcal A_1\mathbin{\dot\cup}
\mathcal A_2\mathbin{\dot\cup}
\mathcal A_3,
$
where
\[
\begin{aligned}
\mathcal A_1
&:=\mathcal D\cap\{|X-Y|<1/4\},\\
\mathcal A_2
&:=\mathcal D\cap\{|X-Y|\ge1/4,\ |Y|\le2\},\\
\mathcal A_3
&:=\mathcal D\cap\{|Y|>2\}.
\end{aligned}
\]
On \(\mathcal A_1\), we have
$
\frac34<|X|<1$
and 
$1<|Y|<\frac54.
$
Thus the integral over
\(\mathcal A_1\) is finite since \(|X-Y|^{1-n}\) is locally
integrable.
On \(\mathcal A_2\), both \(|X-Y|^{1-n}\) and \(|Y|^{-b}\) are
bounded. Hence the corresponding integral is finite because
\(a<n\).
Finally, on \(\mathcal A_3\), we have 
$
|X-Y|\ge |Y|-|X|\ge\frac12|Y|,
$
and therefore
\[
\begin{aligned}
\iint_{\mathcal A_3}
|X|^{-a}|Y|^{-b}|X-Y|^{1-n}\,\dd Y\dd X
\le
C\left(\int_{B_1}|X|^{-a}\,\dd X\right)
 \left(\int_2^\infty r^{-b}\,\dd r\right)
<\infty.
\end{aligned}
\]
Interchanging \(a\) and \(b\), we also have  \(J(b,a)<\infty\), and thus
\begin{equation}\label{tr}
  T_R\le CR^{-1}.
\end{equation}

Next we estimate the boundary term. By
Lemma~\ref{critical-estimates} and the critical identities
\[
\alpha(q+1)=n,
\qquad
\beta(p+1)=n,
\]
we have
\[
u(x)^{q+1}+v(x)^{p+1}\le C|x|^{-n}.
\]
Therefore,
\[
\begin{aligned}
\int_{\partial B_R}
\bigl(u^{q+1}+v^{p+1}\bigr)\,\dd S
\le \frac{C}{R},
\end{aligned}
\]
together with \eqref{tr}, which implies \eqref{cutofferror}.

\medskip\noindent\emph{Step 2: One-dimensional level-set decomposition.}
Fix \(e\in\Sph\). Every point $x\in\R^n$ can be
written as
$
 x=z+te,
$ where $ z\in e^\perp$ and 
$ t\in\R$.
If \(z\ne0\), it follows by \eqref{sourcepointwise} that 
\[
 \begin{aligned}
 f(z+te)&\le C(|z|^2+t^2)^{-(\beta+2)/2},\\
 g(z+te)&\le C(|z|^2+t^2)^{-(\alpha+2)/2}.
 \end{aligned}
\]
Thus the restrictions of \(f\) and \(g\) to \(z+\R e\) are
integrable and tend to zero as \(t\to\pm\infty\). The exceptional
line \(\R e\), corresponding to \(z=0\), is irrelevant because
\(\{0\}\subset e^\perp\) has measure zero.

For \(z\ne0\) and \(h>0\), let \(\mathcal C_f(z,h)\) denote the
family of connected components of
\[
 \{t\in\R:f(z+te)>h\}.
\]
By \eqref{sourcepointwise}, we have   \(f(z+te)\to0\) as \(t\to\pm\infty\), then there exists
\(T=T(z,h)>0\) such that
$
 \{t\in\R:f(z+te)>h\}\subset(-T,T).
$
Thus every \(I\in\mathcal C_f(z,h)\) is a bounded open interval.
Similarly, for \(w\ne0\) and \(\ell>0\), let
\(\mathcal C_g(w,\ell)\) denote the connected components of
\[
 \{\tau\in\R:g(w+\tau e)>\ell\}.
\]
Every \(J\in\mathcal C_g(w,\ell)\) is likewise a bounded open
interval.

Since the components in each of the families \(\mathcal C_f(z,h)\) and \(J\in\mathcal C_g(w,\ell)\)   are pairwise disjoint, for almost every \(z,w\in e^\perp\) and every \(t,\tau\in\R\), we have 
\begin{equation}\label{fg1}
\begin{aligned}
 f(z+te)
 &=\int_0^\infty
 \one_{\{f(z+te)>h\}}\,\dd h
 =\int_0^\infty
   \sum_{I\in\mathcal C_f(z,h)}\one_I(t)\,\dd h,\\
 g(w+\tau e)
&=\int_0^\infty
  \one_{\{g(w+\tau e)>\ell\}}\,\dd\ell
=\int_0^\infty
  \sum_{J\in\mathcal C_g(w,\ell)}
  \one_J(\tau)\,\dd\ell.
\end{aligned}  
\end{equation}

Fix \(a_0>0\) and define
\[
\Phi_{a_0}(c):=
\begin{cases}
0,
 & |c|\le a_0,\\[2mm]
\dfrac{2}{\pi}\operatorname{sgn}(c)
\arctan\bigl(|c|-a_0\bigr),
 & |c|>a_0.
\end{cases}
\]
The function \(\Phi_{a_0}\) is bounded, odd, and nondecreasing. 
For \(I\in\mathcal C_f(z,h)\) and
\(J\in\mathcal C_g(w,\ell)\), define the 
intervals
\[
I^R
:=
I\cap\{t\in\R:z+te\in B_R\},
\qquad
J^R
:=
J\cap\{\tau\in\R:w+\tau e\in B_R\}.
\]
When $I^R$ and $J^R$ are nonempty, denote their midpoints by
\(c_I^R\) and \(c_J^R\), respectively, with empty intersections contributing zero.

We introduce the center-weighted truncated layer functions
\begin{equation}\label{fg2}
  \begin{aligned}
\mathcal J_{f,R}(z+te)
&:=
\int_0^\infty
\sum_{I\in\mathcal C_f(z,h)}
\Phi_{a_0}(c_I^R)\one_{I^R}(t)\,\dd h,\\
\mathcal J_{g,R}(w+\tau e)
&:=
\int_0^\infty
\sum_{J\in\mathcal C_g(w,\ell)}
\Phi_{a_0}(c_J^R)\one_{J^R}(\tau)\,\dd\ell.
\end{aligned}
\end{equation}
Both $\mathcal J_{f,R}(x)$ and $\mathcal J_{g,R}(y)$ are bounded
and compactly supported in $\R^n$. In fact, their support is contained in $B_R$, so it remains to establish boundedness near the origin.
Fix $0<\varepsilon<\min\{a_0,R\}$ and let 
\[
h> M_\varepsilon:=\max_{\partial B_\varepsilon}(f+g).
\]
Consider  an interval
$I\in\mathcal C_f(z,h)$ such that
$
(z+Ie)\cap B_\varepsilon\ne\varnothing.
$
We claim that $z+Ie\subset B_\varepsilon$. Otherwise, there exists $t_{0}$ such that $z+t_{0}e\in \partial  B_\varepsilon$ and $f(z+t_{0}e)>h$,
a contradiction. It follows that 
$
I\subset(-\varepsilon,\varepsilon)$ and
$
|c_I^R|<\varepsilon<a_0,
$
which implies $\Phi_{a_0}(c_I^R)=0$. Consequently, we have 
$
|\mathcal J_{f,R}|\le M_\varepsilon
$ in $B_\varepsilon$, and hence $\mathcal J_{f,R}$ is bounded in $\R^n$. Similarly, we also obtain the boundedness of
$\mathcal J_{g,R}$.

\medskip\noindent\emph{Step 3: Truncated midpoint defect.}
For fixed 
$
z,w\in e^\perp\setminus\{0\}$, $z\ne w$,
 and $ h,\ell>0$, we  choose
$
I\in\mathcal C_f(z,h)$, $
J\in\mathcal C_g(w,\ell).
$
and denote their
midpoints by \(c_I\) and \(c_J\), respectively, so that
\[
I
=
\left(c_I-\frac{|I|}{2},c_I+\frac{|I|}{2}\right),
\qquad
J
=
\left(c_J-\frac{|J|}{2},c_J+\frac{|J|}{2}\right).
\]
Define the function
\[
A_{I,J}(r)
:=
\int_{-|I|/2}^{|I|/2}
\int_{-|J|/2}^{|J|/2}
\Gamma\bigl(z-w+(r+\xi-\eta)e\bigr)
\,\dd\eta\dd\xi.
\]
Since $
K(r):=\Gamma(z-w+re)=\frac{1}{k\om}
\left(|z-w|^{2}+r^2\right)^{-k/2}$,
 the function \(A_{I,J}\) is even and strictly decreasing on
\((0,\infty)\).

For $R>0$, set $t=c_I^R+\xi$ and $
\tau=c_J^R+\eta,$ and  define the truncated interaction integral
$$
\begin{aligned}
P_R:=
\int_{\R^n}
\mathcal J_{f,R}(x)\,
\partial_e I_2(g\one_{B_R})(x)\,\dd x+
\int_{\R^n}
\mathcal J_{g,R}(x)\,
\partial_e I_2(f\one_{B_R})(x)\,\dd x.
\end{aligned}
$$
By \eqref{fg1}, \eqref{fg2}, and Fubini's theorem,
the first integral can be expressed as
\[
\begin{aligned}
&\int_{\R^n}
\mathcal J_{f,R}(x)\,
\partial_e I_2(g\one_{B_R})(x)\,\dd x\\
&=
\int_{e^\perp}\int_{e^\perp}
\int_0^\infty\int_0^\infty
\sum_{I\in\mathcal C_f(z,h)}
\sum_{J\in\mathcal C_g(w,\ell)}
\Phi_{a_0}(c_I^R)\\
&\qquad\times
\left[
\int_{I^R}\int_{J^R}
\partial_e\Gamma
\bigl(z-w+(t-\tau)e\bigr)
\,\dd\tau\dd t
\right]
\,\dd\ell\dd h\dd w\dd z\\
&=
\int_{e^\perp}\int_{e^\perp}
\int_0^\infty\int_0^\infty
\sum_{I\in\mathcal C_f(z,h)}
\sum_{J\in\mathcal C_g(w,\ell)}
\Phi_{a_0}(c_I^R)
A_{I^R,J^R}'(c_I^R-c_J^R)
\,\dd\ell\dd h\dd w\dd z,
\end{aligned}
\]
where the last equality follows from the definition
of $A_{I,J}$.

On the other hand, since 
$
\partial_e\Gamma(y-x)
=-\partial_e\Gamma(x-y),
$
the corresponding part to the
second integral is
$
-\Phi_{a_0}(c_J^R)
A_{I^R,J^R}'(c_I^R-c_J^R).
$
Therefore, we have 
\begin{align}
P_R
={}&
\int_{e^\perp}\int_{e^\perp}
\int_0^\infty\int_0^\infty
\sum_{I\in\mathcal C_f(z,h)}
\sum_{J\in\mathcal C_g(w,\ell)}
\notag\\
&\quad\times
\bigl[
\Phi_{a_0}(c_I^R)
-\Phi_{a_0}(c_J^R)
\bigr]
A_{I^R,J^R}'(c_I^R-c_J^R)
\,\dd\ell\dd h\dd w\dd z.
\label{truncated-defect}
\end{align}
By the evenness and monotonicity of $A_{I^R,J^R}$,
we have
$
rA_{I^R,J^R}'(r)\le0$ for all 
$r\in\R.
$
Since $\Phi_{a_0}$ is nondecreasing, it follows that
\[
\bigl[\Phi_{a_0}(c_I^R)-\Phi_{a_0}(c_J^R)\bigr]
A_{I^R,J^R}'(c_I^R-c_J^R)
\le0.
\]
Hence, by \eqref{truncated-defect}, we have 
$
P_R\le0.
$

\medskip\noindent\emph{Step 4: Removal of the cutoff.}
Write \(u_e:=\partial_eu=\partial_eI_2g\) and \(v_e:=\partial_ev=\partial_eI_2f\), and define
\[
Q_R
:=
\int_{\R^n}
\mathcal J_{f,R}(x)u_e(x)\,\dd x
+
\int_{\R^n}
\mathcal J_{g,R}(x)v_e(x)\,\dd x.
\]
Then 
\begin{equation}\label{PQ1}
  \begin{aligned}
|Q_R-P_R|
=&
\left|\int_{B_R}\int_{B_R^c}
\left[\mathcal J_{f,R}(x)g(y)+\mathcal J_{g,R}(x)f(y)\right]
\partial_e\Gamma(x-y)\,\dd y\dd x\right|\\
&\leq
C\int_{B_R}\int_{B_R^c}
\Bigl[
f(x)g(y)+g(x)f(y)
\Bigr]
|x-y|^{1-n}\,\dd y\dd x= CT_R,
\end{aligned}
\end{equation}
where we have used 
$
\left|\mathcal J_{f,R}\right|
\le f\one_{B_R}$,
$
\left|\mathcal J_{g,R}\right|
\le g\one_{B_R}$,
and
$
|\partial_e\Gamma(x-y)|
\le C|x-y|^{1-n}.
$

We next establish the estimate
\begin{equation}\label{PQ2}
  |Q_R|
\le C\int_{\partial B_R}
\bigl(u^{q+1}+v^{p+1}\bigr)\,\dd S.
\end{equation}
For each $I\in\mathcal C_f(z,h)$ with
$I^R\ne\varnothing$, write
$
I^R=(a_I^R,b_I^R).
$
By the definition of $\mathcal J_{f,R}$ and
Fubini's theorem, we obtain
\[
\begin{aligned}
\int_{\R^n}\mathcal J_{f,R}u_e\,\dd x
=&
\int_{e^\perp}\int_0^\infty
\sum_{I\in\mathcal C_f(z,h)}
\Phi_{a_0}(c_I^R)
\bigl[
u(z+b_I^Re)-u(z+a_I^Re)
\bigr]
\,\dd h\dd z.
\end{aligned}
\]
If $I^R=I$,  we have \(u(z+ae)=u(z+be)=h^{1/q}.\) If $I^R\neq I$, then  $z+a_I^Re\in\partial B_R$, or $z+b_I^Re\in\partial B_{R}$.
For fixed \(z\in e^\perp\), let
$
S_R(z):=(z+\R e)\cap\partial B_R.
$
For \(x\in S_R(z)\),we have
$
0<h<u(x)^q$.
 Thus, using
\(|\Phi_{a_0}|\le1\) and the triangle inequality, we obtain
\[
\begin{aligned}
|\int_{\R^n}\mathcal J_{f,R}u_e\,\dd x|
&\le
\int_{|z|<R}
\sum_{x\in S_R(z)}
\int_0^{u(x)^q}
\bigl[u(x)-h^{1/q}\bigr]
\,\dd h\dd z\\
&=
\frac{1}{q+1}
\int_{|z|<R}
\sum_{x\in S_R(z)}
u(x)^{q+1}\,\dd z\\
&=
\frac{1}{q+1}
\int_{\partial B_R}
u^{q+1}|\nu\cdot e|\,\dd S.
\end{aligned}
\]
Similarly, with \((f,u,q)\) replaced by \((g,v,p)\), we have 
\[
\left|
\int_{\R^n}v_e(x)\mathcal J_{g,R}(x)\,\dd x
\right|
\le
\frac{1}{p+1}
\int_{\partial B_R}
v^{p+1}|\nu\cdot e|\,\dd S.
\]
Therefore, we conclude that 
\[
\begin{aligned}
|Q_R|
&\le
\left|
\int_{\R^n}u_e\mathcal J_{f,R}\,\dd x
\right|
+
\left|
\int_{\R^n}v_e\mathcal J_{g,R}\,\dd x
\right|
\le
C\int_{\partial B_R}
\bigl(u^{q+1}+v^{p+1}\bigr)\,\dd S.
\end{aligned}
\]
Combining \eqref{cutofferror}, \eqref{PQ1} with  \eqref{PQ2}, we
yield
\[
\begin{aligned}
|P_R|
\le |P_R-Q_R|+|Q_R|
\le
CT_R+
C\int_{\partial B_R}
\bigl(u^{q+1}+v^{p+1}\bigr)\,\dd S
\le\frac{C}{R},
\end{aligned}
\]
which implies $P_R\rightarrow0$ as 
$R\to\infty.$

For fixed
$
I\in\mathcal C_f(z,h)$,
and $
J\in\mathcal C_g(w,\ell),
$ define
\[
D_R(I,J)
:=
-\bigl[
\Phi_{a_0}(c_I^R)-\Phi_{a_0}(c_J^R)
\bigr]
A_{I^R,J^R}'(c_I^R-c_J^R),
\]
with \(D_R(I,J):=0\) if either \(I^R\) or \(J^R\) is empty. By
Step~3, we know that 
\[
D_R(I,J)\ge0,
\qquad
-P_R
=
\int_{e^\perp}\int_{e^\perp}
\int_0^\infty\int_0^\infty
\sum_I\sum_J
D_R(I,J)
\,\dd\ell\dd h\dd w\dd z.
\]
 For sufficiently
large \(R\), we have 
\[
I^R=I,
\qquad
J^R=J,
\qquad
c_I^R=c_I,
\qquad
c_J^R=c_J.
\]
It follows that
$
D_R(I,J)\longrightarrow D(I,J),
$ as $R\rightarrow\infty$,
where
$$
\begin{aligned}
D(I,J)
:={}&
-\bigl[
\Phi_{a_0}(c_I)-\Phi_{a_0}(c_J)
\bigr]
A_{I,J}'(c_I-c_J)
\ge0.
\end{aligned}
$$
By Fatou's lemma, we obtain
\[
\begin{aligned}
0
\le
\int_{e^\perp}\int_{e^\perp}
\int_0^\infty\int_0^\infty
\sum_{I\in\mathcal C_f(z,h)}
\sum_{J\in\mathcal C_g(w,\ell)}
D(I,J)
\,\dd\ell\dd h\dd w\dd z
\le
\liminf_{R\to\infty}(-P_R)
=0,
\end{aligned}
\]
and hence $D(I,J)=0$ for almost every $
I\in\mathcal C_f(z,h)$,
and $
J\in\mathcal C_g(w,\ell)$.

\medskip\noindent
\emph{Step 5: Identification of the center.}
Fix \(a_0>0\). Recall from Step~4 that
\begin{equation}\label{midpoint-defect-zero}
\begin{aligned}
0
={}&
\int_{e^\perp}\int_{e^\perp}
\int_0^\infty\int_0^\infty
\sum_{I\in\mathcal C_f(z,h)}
\sum_{J\in\mathcal C_g(w,\ell)}
D(I,J)
\,\dd\ell\dd h\dd w\dd z.
\end{aligned}
\end{equation}

We prove that
\begin{equation}\label{f-midpoint-zero-a0}
\int_{e^\perp}\int_0^\infty
\sum_{I\in\mathcal C_f(z,h)}
\one_{\{|c_I|>a_0\}}
\,\dd h\dd z
=0.
\end{equation}
Suppose, to the contrary, that
\begin{equation}\label{positive-f-midpoints}
\int_{e^\perp}\int_0^\infty
\sum_{I\in\mathcal C_f(z,h)}
\one_{\{|c_I|>a_0\}}
\,\dd h\dd z
>0.
\end{equation}
Choose \(0<\varepsilon<a_0\) and set $
Q_\varepsilon
:=
\max_{\partial B_\varepsilon}g.
$
Since the origin is nonremovable, by Lemma~\ref{Newton}, we may therefore choose
$
x_\varepsilon\in B_\varepsilon\setminus\{0\}$ such that 
$
g(x_\varepsilon)>Q_\varepsilon+1.
$
By continuity, there exists \(\delta>0\) such that
\begin{equation}\label{anchor-ball}
\overline{B_\delta(x_\varepsilon)}
\subset B_\varepsilon\setminus\{0\},
\qquad
g>Q_\varepsilon+1
\quad\text{in }B_\delta(x_\varepsilon).
\end{equation}

Define
\[
W_\varepsilon
:=
\left\{
w\in e^\perp:
(w+\mathbb Re)\cap B_\delta(x_\varepsilon)\ne\varnothing
\right\}.
\]
Then \(W_\varepsilon\) is an open \((n-1)\)-dimensional ball in
\(e^\perp\), and therefore \(|W_\varepsilon|>0\).
For fixed 
$
w\in W_\varepsilon$ and 
$\ell\in(Q_\varepsilon,Q_\varepsilon+1)$,
by Eq.~\eqref{anchor-ball}, then we have 
$$
\{\tau\in\mathbb R:w+\tau e\in B_\delta(x_\varepsilon)\}
\subset
\{\tau\in\mathbb R:g(w+\tau e)>\ell\}.$$
Let \(J\in\mathcal C_g(w,\ell)\) be the connected component containing
$
\{\tau\in\mathbb R:w+\tau e\in B_\delta(x_\varepsilon)\}$.
Thus, 
$w+\tau e\in B_\varepsilon$ for every $\tau\in J$. Otherwise,
there exists
\(\tau_0\) such that 
$
w+\tau_0e\in\partial B_\varepsilon, 
$
which is impossible since 
$
g(w+\tau_0e)
\le Q_\varepsilon
<\ell.
$
Consequently, we have 
$
J\subset(-\varepsilon,\varepsilon)
$ and $
|c_J|<\varepsilon<a_0$,
and hence
$
\Phi_{a_0}(c_J)=0.
$
It follows that
\begin{equation}\label{positive-centered-g-layers}
\begin{aligned}
\int_{e^\perp}\int_0^\infty
\sum_{J\in\mathcal C_g(w,\ell)}
\one_{\{|c_J|<\varepsilon\}}
\,\dd\ell\dd w                                      
\ge
\int_{W_\varepsilon}
\int_{Q_\varepsilon}^{Q_\varepsilon+1}
1\,\dd\ell\dd w
=
|W_\varepsilon|
>0.
\end{aligned}
\end{equation}

For  \(z\ne w\), consider such $
I\in\mathcal C_f(z,h)$ and 
$J\in\mathcal C_g(w,\ell)$ such that 
$|c_I|>a_0$
and $
|c_J|<\varepsilon<a_0$.  Thus, we have 
$  
\begin{aligned}
D(I,J)=-\Phi_{a_0}(c_I)A_{I,J}'(c_I-c_J)
>0.
\end{aligned}
$
Since the diagonal
\[
\{(z,w)\in e^\perp\times e^\perp:z=w\}
\]
has zero product Lebesgue measure, combining
\eqref{positive-f-midpoints} with 
\eqref{positive-centered-g-layers}, we obtain that 
\begin{align*}
0
<&
\int_{e^\perp}\int_{e^\perp}
\int_0^\infty\int_0^\infty
\sum_{I\in\mathcal C_f(z,h)}
\sum_{J\in\mathcal C_g(w,\ell)}
\one_{\{|c_I|>a_0\}}
\one_{\{|c_J|<\varepsilon\}}
D(I,J)
\,\dd\ell\dd h\dd w\dd z                                      \\
\leq&
\int_{e^\perp}\int_{e^\perp}
\int_0^\infty\int_0^\infty
\sum_{I\in\mathcal C_f(z,h)}
\sum_{J\in\mathcal C_g(w,\ell)}
D(I,J)
\,\dd\ell\dd h\dd w\dd z
=0,
\end{align*}
a contradiction.
Hence
Eq.~\eqref{f-midpoint-zero-a0} holds.

Since \(a_0>0\) is arbitrary, and 
$
\{I:c_I\ne0\}
=
\bigcup_{m=1}^\infty
\left\{I:|c_I|>\frac1m\right\}$, $ m\in\mathbb N,
$
we deduce from \eqref{f-midpoint-zero-a0} that
\begin{equation}\label{f-midpoints-zero}
\int_{e^\perp}\int_0^\infty
\sum_{I\in\mathcal C_f(z,h)}
\one_{\{c_I\ne0\}}
\,\dd h\dd z
=0.
\end{equation}
Interchanging \(f\) and \(g\), and using that \(u\) is also unbounded
in every punctured neighborhood of the origin, we obtain
\begin{equation}\label{g-midpoints-zero}
\int_{e^\perp}\int_0^\infty
\sum_{J\in\mathcal C_g(w,\ell)}
\one_{\{c_J\ne0\}}
\,\dd\ell\dd w
=0.
\end{equation}
Thus, by \eqref{f-midpoints-zero}, for almost every
\((z,h)\in e^\perp\times(0,\infty)\), 
we have 
$
\sum_{I\in\mathcal C_f(z,h)}
\one_{\{c_I\ne0\}}
=0$,
which implies
\[
c_I=0
\qquad\text{for every }I\in\mathcal C_f(z,h).
\]
Consequently,  the set
$
\left\{
t\in\mathbb R:f(z+te)>h
\right\}
$
is either empty or a single interval centered at the origin. In particular, there holds 
\[
\one_{\{f(z+te)>h\}}
=
\one_{\{f(z-te)>h\}}
\]
for every \(t\in\mathbb R\) and almost every
\((z,h)\in e^\perp\times(0,\infty)\), then  we have 
\begin{equation}\label{fzt}
\begin{aligned}
f(z+te)
=
\int_0^\infty
\one_{\{f(z+te)>h\}}
\,\dd h
=
\int_0^\infty
\one_{\{f(z-te)>h\}}
\,\dd h
=
f(z-te),
\end{aligned}
\end{equation}
for almost every \((z,t)\in e^\perp\times\mathbb R\).

Applying the same argument to \(g\) and  using
\eqref{g-midpoints-zero}, we yield
\[
g(z+te)=g(z-te)\qquad \text{for~ a.e.} \  (z,t)\in e^\perp\times\mathbb R.
\]
Therefore, 
  we conclude that
\[
u=u\circ\sigma_e,
\qquad
v=v\circ\sigma_e
\qquad\text{in }\mathbb R^n\setminus\{0\},
\]
since
$
\sigma_e(z+te)=z-te$ for  $z\in e^\perp$. We complete the proof.

\end{proof}

\begin{proof}[Proof of Theorem~\ref{symmetry}]

Lemma~\ref{critical-estimates} proves the pointwise bound~\eqref{pointbound}. Since \(e\in\Sph\) is arbitrary,
Lemma~\ref{centered-reflections} shows that \(u\) and \(v\) are
invariant under every hyperplane reflection through the origin.
Thus,  both \(u\) and \(v\) are radial
about the origin.

It remains to prove that both \(u\) and \(v\) are strictly decreasing in the radial variable. By
Lemma~\ref{Newton},
integrating \eqref{distribution} over \(B_r\)
and using the divergence theorem, for every \(r>0\), we obtain
\[
\begin{aligned}
0<\int_{B_r}v^p\,\dd x
=
-\int_{B_r}\Delta u\,\dd x
=
-\int_{\partial B_r}\partial_\nu u\,\dd S
=
-\omega_{n-1}r^{n-1}u'(r),
\end{aligned}
\]
and, similarly,
\[
\begin{aligned}
0<\int_{B_r}u^q\,\dd x
=
-\int_{B_r}\Delta v\,\dd x
=
-\int_{\partial B_r}\partial_\nu v\,\dd S
=
-\omega_{n-1}r^{n-1}v'(r).
\end{aligned}
\]
 Consequently, we have 
\[
u'(r)<0,\qquad v'(r)<0
\qquad (r>0), 
\]
which completes the proof.
\end{proof}

\section{Planar reduction of the radial flow}\label{sec:planar}
Throughout this section, we assume that \eqref{critical}
holds, without imposing any ordering on $p$ and $q$.
By Theorem~\ref{symmetry}, every positive solution of
\eqref{system}--\eqref{critical} with a nonremovable singularity
at the origin is radially symmetric about the origin.
For any positive solution that is radial about the origin, the
Emden--Fowler transformation
\[
U(t)=r^\alpha u(r),\qquad V(t)=r^\beta v(r),\qquad t=\log r
\]
gives a positive global solution of \eqref{ODE}.
Conversely, each positive global solution of \eqref{ODE} gives
a positive radial solution of \eqref{system} on
$\mathbb R^n\setminus\{0\}$ via the inverse transformation. 

In this section, we establish a planar reduction of the radial flow.
We first obtain uniform bounds and compactness for positive global
solutions of \eqref{ODE}. We then prove uniqueness from the data
$(U,U')$ and nondegeneracy of critical points, and construct the
associated planar vector field, which is locally Lipschitz away
from the origin. Set
$z(t)=(U(t),U'(t),V(t),V'(t))$ and define
\[
\K:=\bigl\{
z(t):t\in\R,\ (U,V)\text{ is a positive global solution of }\eqref{ODE}
\bigr\}\cup\{\mathbf{0}\}.
\]

\subsection{Bounds and compactness}

The following lemmas establish the a priori bounds
and compactness needed for the planar reduction.

\begin{lemma}\label{radialbounds}
	There exists a constant $\kappa=\kappa(n)>0$
	such that every positive global solution $(U,V)$
	of \eqref{ODE} satisfies:
	\begin{enumerate}[label=(\roman*),ref=(\roman*)]
		\item\label{uniformbounds}
		$0<U\le\kappa^{-\alpha/2}$ and
		$0<V\le\kappa^{-\beta/2}$.
		
		\item\label{derivativebounds}
		$-\beta U\le U'\le\alpha U$ and
		$-\alpha V\le V'\le\beta V$.
		
		\item\label{crossbounds}
		$U\ge\kappa V^p$ and $V\ge\kappa U^q$.
	\end{enumerate}
	Moreover, for every positive global solution
	$(U,V)$ of \eqref{ODE},
	\begin{equation}\label{wholeGreen}
		U(t)=\frac1{\alpha+\beta}
		\left(
		\int_{-\infty}^t
		\ee^{-\beta(t-s)}V(s)^p\,\dd s
		+\int_t^\infty
		\ee^{\alpha(t-s)}V(s)^p\,\dd s
		\right).
	\end{equation}
	The corresponding formula for $V$ is obtained
	by interchanging $\alpha$ and $\beta$ and
	replacing $V^p$ with $U^q$.
\end{lemma}

\begin{proof}
	Let $(u,v)$ be the positive radial solution of
	\eqref{system} corresponding to $(U,V)$.
	By Theorem~\ref{symmetry} in the nonremovable case
	and \cite[Theorems~1.1 and~2.1]{GHPY}
	in the removable case, we have
	\[
	u'(r)<0,\qquad v'(r)<0
	\quad (r>0),
	\qquad
	(u(r),v(r))\to(0,0)
	\quad\text{as }r\to\infty.
	\]
	
	Set $F(r)=-r^{n-1}u'(r)$.
	Then $F'(r)=r^{n-1}v^p(r)>0$.
	Since $F$ is increasing and $v$ is decreasing,
	\[
	u(r)=\int_r^\infty F(s)s^{1-n}\,\dd s
	\ge\frac{F(r)}{n-2}r^{2-n},
	\qquad
	F(r)\ge\frac{1-2^{-n}}n r^n v^p(r).
	\]
	Consequently,
	\[
	u(r)\ge\kappa r^2v^p(r),
	\qquad
	\kappa=\frac{1-2^{-n}}{n(n-2)},
	\]
	which implies $U\ge\kappa V^p$ since $p\beta=\alpha+2$. Similarly, we have $V\ge\kappa U^q$, which proves \ref{crossbounds}. It follows that \ref{uniformbounds}. Moreover, we have
	\[
	0<-\frac{ru'}u\le n-2.
	\]
	Since $U'/U=\alpha+ru'/u$ and
	$\alpha+\beta=n-2$, we obtain
	\[
	-\beta U\le U'\le\alpha U.
	\]
	An analogous argument applied to $V$ yields $-\alpha V \le V' \le \beta V$, which proves \ref{derivativebounds}.
	
	Let $(U,V)$ be an arbitrary positive global solution
	of \eqref{ODE}. By \ref{uniformbounds}, $U,V$ are bounded on $\mathbb R$. Set
	\[
	\mathcal L:=
	-\frac{\dd^2}{\dd t^2}
	-d\frac{\dd}{\dd t}
	+c.
	\]
	The function
	\[
	G(t)=\frac1{\alpha+\beta}
	\begin{cases}
		e^{-\beta t},&t\ge0,\\
		e^{\alpha t},&t<0
	\end{cases}
	\]
	satisfies $\mathcal LG=\delta_0$ in
	$\mathcal D'(\R)$. 
	Since $V$ is bounded, the function $
	\widetilde U=G*V^p$ is well defined and bounded. Moreover, $\mathcal LG=\delta_0$ yields
	\[
	\mathcal L\widetilde U=V^p
	\quad\text{in }\mathcal D'(\R).
	\]
	Consequently, $\widetilde U$ is a bounded solution
	of $\mathcal L\widetilde U=V^p$. 
	
	Set $W=U-\widetilde U$. Then $W$
	satisfies the following equation
	\[
	-W''-dW'+cW=0.
	\]
	A direct calculation shows  
	\[
	W(t)=ae^{\alpha t}+be^{-\beta t},
	\qquad a,b\in\R.
	\]
	Since $U$ and $\widetilde U$ are bounded
	on $\R$, we have $W\in L^\infty(\R)$. Consequently, $a=b=0$, which implies $U=\widetilde U$. The formula for $V$ follows analogously.
\end{proof}

\begin{lemma}\label{Kcompact}
	The set $\K$ is compact in $\R^4$.
\end{lemma}
\begin{proof}
	Let $\{z_j\}$ be a sequence in $\K$. 
	Since $0\in\K$, we may assume without loss of generality that $z_j\neq 0$ for all $j$.
	Then for each $j$, there exist a positive global solution $(U_j,V_j)$ of \eqref{ODE}
	and $t_j\in\R$ such that
	\[
	z_j = (U_j(t_j), U_j'(t_j), V_j(t_j), V_j'(t_j)).
	\]
	Since \eqref{ODE} is autonomous, after a translation in $t$ we may assume that $t_j=0$ for all $j$.
	
	In view of \ref{uniformbounds} and \ref{derivativebounds}, the sequences $\{U_j\}$, $\{V_j\}$, 
	$\{U_j'\}$, and $\{V_j'\}$ are uniformly bounded on $\mathbb{R}$. 
	It then follows from \eqref{ODE} that $\{U_j''\}$ and $\{V_j''\}$ are also uniformly bounded 
	and equicontinuous on compact intervals. 
	By the Arzel\`a--Ascoli theorem and a standard diagonal argument, passing to a subsequence, 
	we obtain
	\[
	(U_j, V_j) \to (U, V) \quad\text{in } C^2_{\mathrm{loc}}(\mathbb{R}),
	\]
	where $(U, V)$ is a nonnegative entire solution of \eqref{ODE}.
	
	If $U(t_0)=0$ for some $t_0\in\mathbb{R}$, then \ref{derivativebounds} and Gronwall's inequality imply $U \equiv 0$ on $\mathbb{R}$, which in turn forces $V \equiv 0$ via \eqref{ODE}. A similar argument applies if $V$ vanishes at some point.
	Therefore, either $(U, V)\equiv (0,0)$ or $(U, V)$ is a positive global solution of \eqref{ODE}. Consequently, $\K$ is compact in $\R^4$.
\end{proof}

\subsection{Uniqueness and planar reduction}

To obtain a planar reduction of the dynamics
of positive global solutions, we first prove the following
uniqueness and nondegeneracy properties.

\begin{lemma}\label{lem:uniqueness}
	A positive global solution of \eqref{ODE} is uniquely
	determined by either pair of data
	$(U(t_0),U'(t_0))$ or $(V(t_0),V'(t_0))$
	at any $t_0\in\R$. Moreover, every critical point
	of either component of a nonconstant positive global
	solution is nondegenerate.
\end{lemma}
\begin{proof}
	For $\delta\in\{d,-d\}$, We first establish a comparison property for the
	cooperative system
	\begin{equation}\label{cooperative}
		\begin{cases}
			h''+\delta h'-ch=a(t)j,\\
			j''-\delta j'-cj=b(t)h,
		\end{cases}
	\end{equation}
	where $a,b\in C(\R)$ are positive. 
	Set
	\[
	Y(t):=(h(t),h'(t),j(t),j'(t)).
	\]
	Fixing $t_0\in\R$, we claim that if
	$Y(t_0)\in [0,\infty)^4\setminus\{\mathbf{0}\}$, then $Y(t)\in(0,\infty)^4$ for $t>t_0$.
	
	\emph{Case 1. $Y(t_0)\in(0,\infty)^4$.}
	Set
	\[
	\mathcal Z
	:=\{t>t_0:h'(t)=0\text{ or }j'(t)=0\}.
	\]
	If $\mathcal Z\ne\emptyset$, let
	\[
	\tau:=\inf\mathcal Z.
	\]
	By continuity and the positivity of $h'(t_0)$
	and $j'(t_0)$, we have $\tau>t_0$,
	$\tau\in\mathcal Z$, and $
	h'(t),j'(t)>0$ for $t_0\le t<\tau$.
	Consequently, $h(\tau),j(\tau)>0$.
	If $h'(\tau)=0$, then $h''(\tau)\le0$,
	which contradicts
	\[
	h''(\tau)=ch(\tau)+a(\tau)j(\tau)>0.
	\]
	Thus we must have $j'(\tau)=0$. Similarly, we obtain
	\[
	0\ge j''(\tau)=cj(\tau)+b(\tau)h(\tau)>0,
	\]
	which is a contradiction. Hence
	$\mathcal Z=\emptyset$, and $Y(t)\in(0,\infty)^4$ for $t\ge t_0$.
	
	\emph{Case 2. $Y(t_0)\in [0,\infty)^4$ with at least one zero component. }For $\varepsilon>0$,
	let $Y_\varepsilon$ be the solution of the
	first-order system associated with
	\eqref{cooperative}, satisfying
	\[
	Y_\varepsilon(t_0)
	=
	\bigl(
	h(t_0)+\varepsilon,
	h'(t_0)+\varepsilon,
	j(t_0)+\varepsilon,
	j'(t_0)+\varepsilon
	\bigr).
	\]
	By the preceding argument, we have
	$Y_\varepsilon(t)\in(0,\infty)^4$ for $t\ge t_0$.
	By continuous dependence on initial values,
	\[
	Y_\varepsilon\rightarrow Y
	\quad\text{locally uniformly on }[t_0,\infty)
	\]
	as $\varepsilon\to0^+$. Hence $
	Y(t)\in [0,\infty)^4$ for $t\ge t_0$.
	
	Multiplying the first equation of \eqref{cooperative}
	by $e^{\delta t}$ and the second by $e^{-\delta t}$,
	and integrating from $t_0$ to $t$, we obtain
	\begin{align}
		h'(t)
		={}&e^{-\delta(t-t_0)}h'(t_0)
		+\int_{t_0}^t e^{-\delta(t-s)}
		\big(ch(s)+a(s)j(s)\big)\,\dd s,
		\label{eq:h-integral}\\
		j'(t)
		={}&e^{\delta(t-t_0)}j'(t_0)
		+\int_{t_0}^t e^{\delta(t-s)}
		\big(cj(s)+b(s)h(s)\big)\,\dd s.
		\label{eq:j-integral}
	\end{align}
	If $(h(t_0),h'(t_0))\ne(0,0)$, the nonnegativity
	of $h$ and $h'$ gives
	\[
	h(t)=h(t_0)+\int_{t_0}^t h'(s)\,\dd s>0,
	\qquad t>t_0.
	\]
	By \eqref{eq:j-integral}, we have
	\[
	j'(t)\ge
	\int_{t_0}^t e^{\delta(t-s)}b(s)h(s)\,\dd s>0,
	\qquad t>t_0.
	\]
	Consequently,
	\[
	j(t)=j(t_0)+\int_{t_0}^t j'(s)\,\dd s>0,
	\qquad t>t_0.
	\]
	By \eqref{eq:h-integral}, we obtain
	\[
	h'(t)\ge
	\int_{t_0}^t e^{-\delta(t-s)}a(s)j(s)\,\dd s>0,
	\qquad t>t_0.
	\]
	Thus $Y(t)\in(0,\infty)^4$ for $t>t_0$. If $h(t_0)=h'(t_0)=0$, then
	$(j(t_0),j'(t_0))\ne(0,0)$. Similarly, we obtain $Y(t)\in(0,\infty)^4$ for $ t>t_0$, which proves the claim.
	
	We next show that every solution of \eqref{cooperative}
	with $Y(t_0)\in[0,\infty)^4\setminus\{\mathbf{0}\}$
	satisfies $j(t)\to+\infty$ as $t\to+\infty$.
	Fixing $t_1>t_0$, we have
	\[
	j(t_1)>0,\quad j'(t_1)>0,\quad
	b(t)h(t)>0\quad\text{ for }t\ge t_1.
	\]
	The characteristic roots of
	\[
	j''-\delta j'-cj=0
	\]
	are $\lambda$ and $-\mu$, where
	\[
	\lambda=\frac{\alpha+\beta+\delta}{2},
	\qquad
	\mu=\frac{\alpha+\beta-\delta}{2}.
	\]
	In particular, $\lambda,\mu>0$ and
	$\lambda+\mu=\alpha+\beta$.
	Then the variation-of-constants formula yields
	\begin{align*}
		j(t)
		&=
		\frac{j'(t_1)+\mu j(t_1)}{\alpha+\beta}
		e^{\lambda(t-t_1)}
		+\frac{\lambda j(t_1)-j'(t_1)}
		{\alpha+\beta}e^{-\mu(t-t_1)}
		\\
		&\quad+\frac1{\alpha+\beta}
		\int_{t_1}^t
		\big(e^{\lambda(t-s)}-e^{-\mu(t-s)}\big)
		b(s)h(s)\,\dd s\\
		&\ge \frac{j'(t_1)+\mu j(t_1)}{\alpha+\beta}
		e^{\lambda(t-t_1)}
		+\frac{\lambda j(t_1)-j'(t_1)}
		{\alpha+\beta}e^{-\mu(t-t_1)}\to\infty\quad\text{ as }t\to\infty.
	\end{align*}
	
	Now we prove that if $Y$ is bounded on $\mathbb R$ and $h(t_0)=h'(t_0)=0$,
	then
	\begin{equation}\label{comparisonconclusion}
		j(t_0)=j'(t_0)=0.
	\end{equation}
	Suppose to the contrary that $(j(t_0),j'(t_0))\ne(0,0)$.
	Choose $\eta,\sigma\in\{-1,1\}$ such that $\eta j(t_0)\ge0$ and $\eta\sigma j'(t_0)\ge0$.
	Define
	\[
	\begin{aligned}
		\widetilde h(t)
		&=\eta h\big(t_0+\sigma(t-t_0)\big),\\
		\widetilde j(t)
		&=\eta j\big(t_0+\sigma(t-t_0)\big).
	\end{aligned}
	\]
	Then $(\widetilde h,\widetilde j)$ satisfies
	\[
	\begin{cases}
		\widetilde h''+\sigma\delta\widetilde h'-c\widetilde h
		=a\big(t_0+\sigma(t-t_0)\big)\widetilde j,\\
		\widetilde j''-\sigma\delta\widetilde j'-c\widetilde j
		=b\big(t_0+\sigma(t-t_0)\big)\widetilde h.
	\end{cases}
	\]
	Moreover,
	\[
	\begin{aligned}
		\big(
		\widetilde h(t_0),\widetilde h'(t_0),
		\widetilde j(t_0),\widetilde j'(t_0)
		\big)=
		\big(0,0,\eta j(t_0),\eta\sigma j'(t_0)\big)
		\in[0,\infty)^4\setminus\{0\}.
	\end{aligned}
	\]
	Since $\sigma\delta\in\{d,-d\}$, we have $\widetilde j(t)\rightarrow+\infty$as $t\to+\infty$.
	This contradicts the boundedness of $Y$ on
	$\mathbb R$, and hence \eqref{comparisonconclusion} holds.
	
	Let $(U_1,V_1)$ and $(U_2,V_2)$ be two positive global
	solutions of \eqref{ODE}. We show that
	$(U_1,V_1)\equiv(U_2,V_2)$ if
	$U_1(t_0)=U_2(t_0)$ and $U_1'(t_0)=U_2'(t_0)$.
	The case $V_1(t_0)=V_2(t_0)$ and
	$V_1'(t_0)=V_2'(t_0)$ is analogous.
	Let
	\[
	h=U_1-U_2,\qquad j=V_2-V_1.
	\]
	Then $h$ and $j$ satisfy \eqref{cooperative} with $\delta=d$ and
	\[
	\begin{aligned}
		a(t)
		&=p\int_0^1
		\big((1-s)V_1(t)+sV_2(t)\big)^{p-1}\,\dd s,\\
		b(t)
		&=q\int_0^1
		\big((1-s)U_2(t)+sU_1(t)\big)^{q-1}\,\dd s.
	\end{aligned}
	\]
	Moreover, $a,b\in C(\R)$ are positive.
	By Lemma~\ref{radialbounds}, we have 
	$h,j,h',j'\in L^\infty(\mathbb R)$. Combining this with $h(t_0)=h'(t_0)=0$, we get $j(t_0)=j'(t_0)=0$. By uniqueness of the initial value problem associated
	with \eqref{ODE}, we obtain $(U_1,V_1)\equiv(U_2,V_2)$.
	
	It remains to prove the nondegeneracy of critical
	points. Let $(U,V)$ be a nonconstant positive global solution of \eqref{ODE}.
	Differentiating \eqref{ODE} and setting
	\[
	\bar h=U',\qquad \bar j=-V',
	\]
	we obtain
	\[
	\begin{cases}
		\bar h''+d\bar h'-c\bar h=pV^{p-1}\bar j,\\
		\bar j''-d\bar j'-c\bar j=qU^{q-1}\bar h.
	\end{cases}
	\]
	By Lemma~\ref{radialbounds} and \eqref{ODE},
	$\bar h,\bar j,\bar h',\bar j'$ are bounded on
	$\mathbb R$. Suppose to the contrary that $U'(t_0)=U''(t_0)=0$ for some $t_0\in\R$.
	Then $\bar h(t_0)=\bar h'(t_0)=0$. It follows that $
	\bar j(t_0)=\bar j'(t_0)=0$, i.e., $
	V'(t_0)=V''(t_0)=0$.
	Then \eqref{ODE} yields
	\[
	V^p(t_0)=cU(t_0),\qquad
	U^q(t_0)=cV(t_0).
	\]
	Hence $
	\big(U(t_0),0,V(t_0),0\big)$ is an equilibrium of the first-order system
	associated with \eqref{ODE}. By uniqueness of the initial value
	problem for this first-order system, $(U,V)$
	must be constant. This contradicts the assumption that $(U,V)$ is
	nonconstant. Hence every critical point of $U$
	is nondegenerate. Similarly, every critical point of $V$ is nondegenerate.
\end{proof}

Define the projection $\pi:\K\to\R^2$ by
$\pi(U,U',V,V')=(U,U')$. We have

\begin{lemma}\label{lem:projection}
	The map $\pi$ is a homeomorphism from $\K$ onto
	$\pi\K$. Moreover, its inverse $\pi^{-1}$ is locally
	Lipschitz on $\pi\K\setminus\{0\}$.
\end{lemma}
\begin{proof}
	By Lemma~\ref{lem:uniqueness}, $\pi$ is injective. Since $\K$ is compact
	and $\pi$ is continuous, $\pi\K$ is compact.
	Thus $\pi:\K\to\pi\K$ is a homeomorphism.
	
	We now prove that $\pi^{-1}$ is locally Lipschitz on
	$\pi\K\setminus\{0\}$. Assume by contradiction that there exists $x_*\in\pi\K\setminus\{0\}$ such that
	$\pi^{-1}$ is not locally Lipschitz at $x_*$.
	Set $z_*=\pi^{-1}(x_*)$.
	Then there exist two sequences $\{z_k\},\{\bar z_k\}\subset\K$ such that
	\[
	z_k,\bar z_k\to z_*,
	\qquad
	\delta_k:=|z_k-\bar z_k|>0,
	\qquad
	\frac{|\pi z_k-\pi\bar z_k|}{\delta_k}
	\rightarrow0.
	\]
	Let $(U_*,V_*)$, $(U_k,V_k)$, and
	$(\overline U_k,\overline V_k)$ denote the positive global
	solutions of \eqref{ODE} passing through
	$z_*$, $z_k$, and $\bar z_k$, respectively,
	at $t=0$.
	By continuous dependence on initial data, we have
	\[
	(U_k,V_k),\;(\overline U_k,\overline V_k)
	\rightarrow (U_*,V_*)
	\quad\text{in }C^2_{\mathrm{loc}}(\R).
	\]
	
	Set
	\[
	H_k=\frac{U_k-\overline U_k}{\delta_k},
	\qquad
	J_k=\frac{\overline V_k-V_k}{\delta_k}.
	\]
	A direct calculation shows that 
	\[
	\begin{cases}
		H_k''+dH_k'-cH_k=a_k(t)J_k,\\
		J_k''-dJ_k'-cJ_k=b_k(t)H_k,
	\end{cases}
	\]
	where
	\[
	\begin{aligned}
		a_k(t)
		&=p\int_0^1
		\big((1-s)V_k(t)+s\overline V_k(t)\big)^{p-1}\,\dd s,\\
		b_k(t)
		&=q\int_0^1
		\big((1-s)\overline U_k(t)+sU_k(t)\big)^{q-1}\,\dd s.
	\end{aligned}
	\]
	In particular, $a_k\to pV_*^{p-1}$ and
	$b_k\to qU_*^{q-1}$ locally uniformly in $\R$.
	By the definition of $\delta_k$, we see
	\[
	H_k^2(0)+H_k'^2(0)+J_k^2(0)+J_k'^2(0)=1.
	\]
	Moreover,
	\[
	H_k^2(0)+H_k'^2(0)=\frac{|\pi z_k-\pi\bar z_k|^2}{\delta_k^2}\rightarrow0.
	\]
	Passing to a subsequence, we may assume that there exist constants $\gamma_0, \gamma_1$ such that $\gamma_0^2+\gamma_1^2=1$ and 
	\[
	(J_k(0),J_k'(0))\rightarrow(\gamma_0,\gamma_1).
	\]
	By the continuous dependence of solutions of linear
	ODE systems on their coefficients and initial data, there exists a pair $(H,J)$ such that
	\[
	(H_k,J_k)\rightarrow(H,J)
	\quad\text{in }C^2_{\mathrm{loc}}(\R),
	\]
	where $(H,J)$ satisfies
	\[
	\begin{cases}
		H''+dH'-cH=pV_*^{p-1}J,\\
		J''-dJ'-cJ=qU_*^{q-1}H,
	\end{cases}
	\]
	with $H(0)=H'(0)=0$ and
	$J(0)=\gamma_0$, $J'(0)=\gamma_1$.
	
	Choose $\varepsilon,\sigma\in\{-1,1\}$ such that $\varepsilon\gamma_0\ge0$ and $\varepsilon\sigma\gamma_1\ge0$.
	Set
	\[
	(\widehat H_k(t),\widehat J_k(t))=(\varepsilon H_k(\sigma t),\varepsilon J_k(\sigma t)),\qquad (\widehat H(t),\widehat J(t))=(\varepsilon H(\sigma t),\varepsilon J(\sigma t)).
	\]
	Then $
	(\widehat H_k,\widehat J_k)
	\rightarrow(\widehat H,\widehat J)$ in $C^2_{\mathrm{loc}}(\R)$, where $(\widehat H,\widehat J)$ satisfies
	\eqref{cooperative} with $\delta=\sigma d$,
	$a(t)=pV_*(\sigma t)^{p-1}$, and
	$b(t)=qU_*(\sigma t)^{q-1}$.
	Moreover,
	\[
	\big(\widehat H(0),\widehat H'(0),
	\widehat J(0),\widehat J'(0)\big)
	=
	(0,0,\varepsilon\gamma_0,\varepsilon\sigma\gamma_1)
	\in[0,\infty)^4\setminus\{0\}.
	\]
	By the proof of Lemma~\ref{lem:uniqueness}, we obtain that 
	$\widehat H,\widehat H',
	\widehat J,\widehat J'$ are positive on $(0,\infty)$.
	Fix $t_1>0$. The local $C^2$ convergence implies that,
	for all sufficiently large $k$,
	\[
	\widehat H_k(t_1),\widehat H_k'(t_1),
	\widehat J_k(t_1),\widehat J_k'(t_1)>0.
	\]
	Since $(\widehat H_k,\widehat J_k)$ also satisfies
	\eqref{cooperative}, with $\delta=\sigma d$ and
	positive coefficients $a_k(\sigma t)$ and $b_k(\sigma t)$, it follows that $
	\widehat J_k(t)\rightarrow+\infty$ as $t\to+\infty$.
	
	On the other hand, Lemma~\ref{radialbounds} yields
	\[
	|\widehat J_k(t)|=\frac{|\varepsilon\big(\overline V_k(\sigma t)-V_k(\sigma t)\big)|}
	{\delta_k}
	\le
	\frac{\|V_k\|_{L^\infty(\R)}
		+\|\overline V_k\|_{L^\infty(\R)}}{\delta_k}
	\le\frac{2\kappa^{-\beta/2}}{\delta_k},
	\qquad t\in\R,
	\]
	which is a contradiction.
	Therefore, $\pi^{-1}$ is locally Lipschitz on
	$\pi\K\setminus\{0\}$.
\end{proof}

For every positive global solution $(U,V)$ of
\eqref{ODE} and every $t\in\R$, define
$\Phi:\pi\K\to\R$ by
\[
\Phi(U(t),U'(t))=V(t).
\]
By Lemma~\ref{lem:uniqueness}, the definition of $\Phi$
is well--defined. Lemma~\ref{lem:projection} yields that $\Phi$ is locally Lipschitz on $\pi\K\setminus\{0\}$. Moreover, $\Phi(0,0)=0$ and $\Phi>0$ on $\pi\K\setminus\{0\}$.
Define
\[
\mathcal F(x,y)
=\big(y,-dy+cx-\Phi^p(x,y)\big),
\qquad (x,y)\in\pi\K\setminus\{0\}.
\]
By the definition of $\Phi$ and \eqref{ODE}, every
positive global solution $(U,V)$ satisfies
\[
\frac{\dd}{\dd t}(U,U')
=\big(U',-dU'+cU-\Phi^p(U,U')\big)
=\mathcal F(U,U').
\]
Since $\Phi$ is positive and locally Lipschitz on
$\pi\K\setminus\{0\}$, the map $\mathcal F$ is
locally Lipschitz on $\pi\K\setminus\{0\}$. We have the following lemma.
\begin{lemma}\label{lem:planar}
	For every nonempty compact set $K_0\subset\pi\mathcal K\setminus\{0\}$,
	there exist a bounded open set $\mathcal O\subset\mathbb R^2$
	satisfying $K_0\subset\mathcal O$ and
	$0\notin\overline{\mathcal O}$, and a globally
	Lipschitz map
	$\widetilde{\mathcal F}:\mathbb R^2\to\mathbb R^2$
	such that
	\begin{equation}\label{eq:Lipschitz}
		\widetilde{\mathcal F}(x,y)=\mathcal F(x,y),
		\qquad (x,y)\in\mathcal O\cap\pi\mathcal K.
	\end{equation}
	Moreover, $(c^{\alpha/2},0)$ is the unique equilibrium
	of $\mathcal F$ in $\pi\K\setminus\{0\}$.
\end{lemma}
\begin{proof}
	Since $K_0$ is compact and $0\notin K_0$,
	we may choose a bounded open neighborhood
	$\mathcal O\subset\R^2$ of $K_0$ such that
	$0\notin\overline{\mathcal O}$. Set $
	E=\pi\K\cap\overline{\mathcal O}$.  Since $\pi\K$ is compact and
	$0\notin\overline{\mathcal O}$, we have
	$E\subset\pi\K\setminus\{0\}$ is compact. Letting
	\[
	\mathcal F_2(x,y)=-dy+cx-\Phi^p(x,y),
	\]
	we see that $\mathcal F_2$ is locally Lipschitz
	on $\pi\K\setminus\{0\}$.
	For each $\xi\in E$, there exist $r_\xi>0$
	and $L_\xi>0$ such that
	\[
	|\mathcal F_2(\eta)-\mathcal F_2(\zeta)|
	\le L_\xi|\eta-\zeta|,\qquad\forall \eta,\zeta\in E\cap B_{r_\xi}(\xi).
	\]
	Since $E$ is compact, there exist
	$\xi_1,\ldots,\xi_N\in E$ such that $
	E\subset\bigcup_{i=1}^N B_{r_{\xi_i}/2}(\xi_i)$.
	By the Lebesgue number lemma, there exists $\rho>0$
	such that, for every $\xi\in E$,
	\[
	E\cap B_\rho(\xi)\subset B_{r_{\xi_i}/2}(\xi_i)
	\quad\text{for some }i\in\{1,\ldots,N\}.
	\]
	Set
	\[
	L_0=\max_{1\le i\le N}L_{\xi_i},\qquad
	M=\max_{\xi\in E}|\mathcal F_2(\xi)|.
	\]
	For any $\eta,\zeta\in E$ with
	$|\eta-\zeta|<\rho$, there exists $i\in\{1,\ldots,N\}$ such that $
	\eta,\zeta\in B_{r_{\xi_i}/2}(\xi_i)
	\subset B_{r_{\xi_i}}(\xi_i)$.
	Hence
	\[
	|\mathcal F_2(\eta)-\mathcal F_2(\zeta)|
	\le L_{\xi_i}|\eta-\zeta|
	\le L_0|\eta-\zeta|.
	\]
	On the other hand, if $|\eta-\zeta|\ge\rho$,
	then
	\[
	|\mathcal F_2(\eta)-\mathcal F_2(\zeta)|
	\le 2M
	\le \frac{2M}{\rho}|\eta-\zeta|.
	\]
	Therefore, $\mathcal F_2$ is Lipschitz on $E$
	with Lipschitz constant $
	L_E=\max\left\{L_0,2M/\rho\right\}$. 
	Define
	\[
	\widetilde{\mathcal F}(x,y)
	=\big(y,\widetilde{\mathcal F}_2(x,y)\big),
	\]
	where 
	\[
	\widetilde{\mathcal F}_2(\xi)
	:=\inf_{\eta\in E}
	\left\{
	\mathcal F_2(\eta)+L_E|\xi-\eta|
	\right\},
	\qquad \xi\in\R^2.
	\]
	It is straightforward to verify that $\widetilde{\mathcal F}_2$ is globally $L_E$-Lipschitz on $\mathbb{R}^2$, and so $\widetilde{\mathcal F}:\R^2\to\R^2$
	is globally Lipschitz. Since $\mathcal O\cap\pi\K\subset E$ and
	$\widetilde{\mathcal F}_2=\mathcal F_2$ on $E$,
	we obtain \eqref{eq:Lipschitz}.

	We now show that $(c^{\alpha/2},0)$ is an equilibrium
	of $\mathcal F$. The constant pair
	$(U,V)=(c^{\alpha/2},c^{\beta/2})$
	is a positive global solution of \eqref{ODE}. Hence $(c^{\alpha/2},0,c^{\beta/2},0)\in\K$. By the definitions of $\pi$ and $\Phi$, we obtain $(c^{\alpha/2},0)\in\pi\K$ and $\Phi(c^{\alpha/2},0)=c^{\beta/2}$. Thus $\mathcal F(c^{\alpha/2},0)=0$.
	
	Let $(\xi,0)\in\pi\K\setminus\{0\}$ be an
	equilibrium of $\mathcal F$. Then $c\xi=\Phi^p(\xi,0)$.
	By the definition of
	$\K$ and $\pi$, there exists a positive global
	solution $(U,V)$ of \eqref{ODE} such that $(U(t_0),U'(t_0))=(\xi,0)$ for some $t_0\in\R$.
	By the definition of $\Phi$, we have
	\[
	\Phi(\xi,0)=\Phi(U(t_0),U'(t_0))=V(t_0).
	\]
	It follows from \eqref{ODE} that
	\[
	U''(t_0)=c\xi-V^p(t_0)
	=c\xi-\Phi^p(\xi,0)=0.
	\]
	Thus $U'(t_0)=U''(t_0)=0$. By
	Lemma~\ref{lem:uniqueness}, $(U,V)$ must be constant. A direct calculation shows $\xi=c^{\alpha/2}$.
\end{proof}

\section{The global Delaunay family}\label{sec:Delaunay}
We establish a global family of positive periodic solutions of
\eqref{ODE} and characterize the asymptotic behavior of the remaining positive trajectories in $\K$. We begin by constructing periodic orbits near the positive equilibrium $(A,0,B,0)$ and proving their local continuation around arbitrary nonconstant periodic orbits. The projections of these families provide open disk and annular neighborhoods in the $(U,U')$ plane, leading to a periodicity criterion based on limit sets. We then obtain a global parametrization of the periodic orbits by the minimum of $U$ and determine the corresponding energy range and limiting periods. Finally, we prove that every positive trajectory which is neither constant nor periodic converges to zero as $t\to\pm\infty$.

\subsection{Local families of positive periodic orbits}
\label{sec:localperiodic}

We first construct two types of local families of positive periodic orbits. By the Lyapunov center theorem
\cite[Section~I.11.1]{Kielhofer},we  obtain small periodic orbits near the positive equilibrium. Using the Poincar\'e continuation theorem\cite[Chapter~4, Section~4.1, Proposition~2]{HZ}, we also prove that every positive nonconstant periodic
orbit belongs to a smooth local family parametrized by the conserved
energy.

\Needspace{7\baselineskip}
\begin{proposition}
	\label{localperiodic}
	Let $\pi$ be the projection defined in
	Lemma~\ref{lem:projection}.
	The following assertions hold.
	
	\begin{enumerate}[label=\textup{(\roman*)}]
		\item 
		There exist $a_*>0$, $\eta>0$, and a smooth family of
		positive nonconstant periodic solutions
		\[
		z_a=(U_a,U_a',V_a,V_a'),
		\qquad 0<a<a_*,
		\]
		with orbits $\gamma_a:=z_a(\R)$ and least periods $T(a)$,
		such that
		\[
		\sup_{t\in\R}
		\bigl|z_a(t)-(A,0,B,0)\bigr|\longrightarrow0,
		\qquad
		T(a)\longrightarrow\frac{2\pi}{\omega}
		\qquad\text{as }a\rightarrow0.
		\]
		The set
		$
		\{(A,0)\}\cup\bigcup_{0<a<a_*}\pi\gamma_a
		$
		is an open neighborhood of $(A,0)$, and
		\[
		\left\{\min_{t\in\R}U_a(t):0<a<a_*\right\}
		=(A-\eta,A).
		\]
		Each of $U_a$ and $V_a$ has exactly one nondegenerate
		maximum and one nondegenerate minimum in each least period.
		
		\item 
		Let $\gamma_0\subset\K$ be a nonconstant periodic orbit
		with energy $E_0$.
		There exist $\delta>0$ and a smooth family of positive
		nonconstant periodic orbits
		$
		\{\gamma_E:|E-E_0|<\delta\}
		$
		such that
		\[
		\gamma_{E_0}=\gamma_0,
		\qquad
		\cH(z)=E
		\quad\text{for every }z\in\gamma_E.
		\]
		Moreover,
		$
		\bigcup_{|E-E_0|<\delta}\pi\gamma_E
		$
		is an open annulus containing $\pi\gamma_0$.
	\end{enumerate}
\end{proposition}

\begin{proof}
	Consider the Lagrangian
	\begin{equation}\label{lagrangian}
		L=U'V'+\frac d2(UV'-U'V)+cUV
		-\frac{U^{q+1}}{q+1}-\frac{V^{p+1}}{p+1}.
	\end{equation}
	Its Euler--Lagrange equations coincide with
	\eqref{ODE}. The associated conjugate momenta
	and velocity Hessian are
	given by
	\[
	P_U=V'-\frac d2V,\qquad
	P_V=U'+\frac d2U,\qquad
	D^2_{(U',V')}L=
	\begin{pmatrix}
		0&1\\
		1&0
	\end{pmatrix}.
	\]
	Thus,
	the Legendre transform is invertible, and the
	corresponding Hamiltonian is
	\[
	P_UU'+P_VV'-L=\cH.
	\]
	In the canonical coordinates $(U,V,P_U,P_V)$,
	the system \eqref{ODE} admits a canonical
	Hamiltonian formulation, and its linearized
	flow is symplectic.
	
	\medskip\noindent\emph{Proof of \textup{(i)}: Small periodic orbits.}
	At the equilibrium $(A,0,B,0)$, the characteristic
	equation of the linearized system is
	\begin{equation}\label{characteristic}
		\lambda^4-(2c+d^2)\lambda^2-(pq-1)c^2=0.
	\end{equation}
	Since $pq>1$, its roots are simple and given by
	\[
	\pm\rho,\qquad \pm i\omega,
	\qquad \rho,\omega>0,
	\]
	where $\omega$ is defined in \eqref{frequency}.
	In particular, no integer multiple $mi\omega$
	with $|m|\ge2$ is an eigenvalue.
	Applying the Lyapunov center
	theorem \cite[Section~I.11.1]{Kielhofer},
	there exist $a_*>0$ and a smooth
	family of nonconstant periodic solutions
	\[
	z_a=(U_a,U_a',V_a,V_a'),
	\qquad 0<a<a_*,
	\]
	whose periods satisfy
	\[
	T(a)\longrightarrow\frac{2\pi}{\omega}
	\qquad\text{as }a\rightarrow0.
	\]

	Let $(\xi,\eta)$ be the position components of an
	eigenvector corresponding to $i\omega$, then
	\[
	(-\omega^2+i d\omega-c)\xi+pB^{p-1}\eta=0,
	\qquad
	qA^{q-1}\xi+(-\omega^2-i d\omega-c)\eta=0.
	\]
	Since $pB^{p-1}>0$ and $qA^{q-1}>0$, either $\xi=0$
	or $\eta=0$ would imply $\xi=\eta=0$, contradicting the
	nontriviality of the eigenvector. Hence
	$
	\xi\eta\ne0.
	$
	After a suitable normalization of the amplitude
	and a choice of phase, the rescaled profiles
	\[
	\widehat U_a(\theta)
	:=U_a\left(\frac{T(a)}{2\pi}\theta\right),
	\qquad
	\widehat V_a(\theta)
	:=V_a\left(\frac{T(a)}{2\pi}\theta\right)
	\]
	extend smoothly to $a=0$ and admit the expansions
	\begin{equation}\label{uva}
		\begin{aligned}
			\widehat U_a(\theta)
			&=A+a\cos\theta+O(a^2),\\
			\widehat V_a(\theta)
			&=B+ab_*\cos(\theta+\vartheta_*)+O(a^2),
		\end{aligned}	
	\end{equation}
	in $C^2(\mathbb R/(2\pi\mathbb Z))$,
	where $b_*>0$ and $\vartheta_*\in\mathbb R$
	are constants.
	It follows from \eqref{uva} that for all
	sufficiently small $a>0$,
	$
	U_a(t)\ge\frac A2$ and 
	$V_a(t)\ge\frac B2$ for all
	$t\in\mathbb R$.
	Consequently, $\gamma_a:=z_a(\mathbb R)$
	is a positive periodic orbit in $\K$.
	
	Differentiating \eqref{uva}, we obtain
	\[
	a^{-1}\partial_\theta\widehat U_a
	\longrightarrow -\sin\theta,
	\qquad
	a^{-1}\partial_\theta\widehat V_a
	\longrightarrow
	-b_*\sin(\theta+\vartheta_*)
	\]
	in $C^1(\mathbb R/(2\pi\mathbb Z))$ as $a\rightarrow0$.
	Since both $-\sin\theta$ and
	$-b_*\sin(\theta+\vartheta_*)$ have exactly two
	simple zeros on $\mathbb R/(2\pi\mathbb Z)$, for all sufficiently small $a>0$,
	both $\partial_\theta\widehat U_a$ and
	$\partial_\theta\widehat V_a$ have exactly two
	simple zeros in one period.
	Consequently, each of $U_a$ and $V_a$ has exactly
	one nondegenerate maximum and one nondegenerate
	minimum in $[0,T(a))$.
	
	We claim that $T(a)$ is the least period of both
	components. Indeed, let $P_a^U$ denote the least
	positive period of $U_a$. Since $T(a)$ is also a
	period, there exists an integer $m\ge1$ such that
	$
	T(a)=mP_a^U.
	$
	If $m\ge2$, the critical points of $U_a$ would
	repeat at least twice in $[0,T(a))$, contradicting
	the fact that $U_a$ has exactly two critical
	points there. Hence $P_a^U=T(a)$.
	The same argument yields $P_a^V=T(a)$.
	Therefore, $T(a)$ is the least positive period
	of both $U_a$ and $V_a$, and hence of the full
	periodic orbit $z_a$.

	The smooth extension of
	$a^{-1}\partial_\theta\widehat U_a$
	has a simple zero at $(a,\theta)=(0,\pi)$.
	By the implicit function theorem, there exists
	a smooth function $\theta_-(a)$, with
	$\theta_-(0)=\pi$, such that
	$
	\partial_\theta\widehat U_a(\theta_-(a))=0.
	$
	For sufficiently small $a>0$, this critical point
	is the unique minimum of $\widehat U_a$.
	Consequently,
	\[
	m(a):=\widehat U_a(\theta_-(a))
	=A-a+O(a^2),
	\qquad
	m'(a)=-1+O(a)<0,
	\]
	if $a_*$ small enough.
	Since $m(a)\to A$ as $a\rightarrow0$,
	there exists $\eta>0$ such that
	\[
	m((0,a_*))=(A-\eta,A).
	\]
	In particular, distinct values of $a$
	correspond to distinct periodic orbits.
	
	Let
	$
	\mathcal D
	:=\{(x,y)\in\mathbb R^2:x^2+y^2<a_*^2\},
	$
	and define $\mathcal F:\mathcal D\to\mathbb R^2$ by
	\[
	\mathcal F(a\cos\theta,a\sin\theta)
	:=
	\begin{cases}
		\displaystyle
		\pi z_a\left(\frac{T(a)}{2\pi}\theta\right),
		&0<a<a_*,\\[2mm]
		(A,0),&a=0.
	\end{cases}
	\]
	This map is well defined by the periodicity of $z_a$.
	Moreover, $\mathcal F$ is continuous away from
	the origin by smooth dependence, and continuous
	at the origin by the uniform convergence
	$z_a\to(A,0,B,0)$ as $a\rightarrow0$.

	We now prove that $\mathcal F$ is injective. Let
	$
	x=a(\cos\theta,\sin\theta)$
	and $
	\widetilde x=\widetilde a
	(\cos\widetilde\theta,\sin\widetilde\theta)
	$
	be points of $\mathcal D$ such that
	$
	\mathcal F(x)=\mathcal F(\widetilde x).
	$
	Assume first that $a,\widetilde a>0$. By the definition of
	$\mathcal F$, we have
	$
	\pi z_a\left(\frac{T(a)}{2\pi}\theta\right)
	=
	\pi z_{\widetilde a}
	\left(\frac{T(\widetilde a)}{2\pi}\widetilde\theta\right).
	$
	Since $\pi$ is injective on $\K$ by
	Lemma~\ref{lem:projection}, then
	$
	z_a\left(\frac{T(a)}{2\pi}\theta\right)
	=
	z_{\widetilde a}
	\left(\frac{T(\widetilde a)}{2\pi}\widetilde\theta\right).
	$
	By the uniqueness for the initial-value problem associated with
	\eqref{ODE}, we obtain  
	$
	\gamma_a=\gamma_{\widetilde a},
	$
	and hence
	$
	m(a)=m(\widetilde a).
	$
	By the strict monotonicity of $m$, we obtain
	$a=\widetilde a$. Consequently,
	\[
	z_a\left(\frac{T(a)}{2\pi}\theta\right)
	=
	z_a\left(\frac{T(a)}{2\pi}\widetilde\theta\right).
	\]
	Since  $T(a)$ is the least period of $z_a$, we obtain
	$
	\theta-\widetilde\theta\in2\pi\mathbb Z.
	$
	Therefore \(x=\widetilde x\).
	
	It remains to consider the origin. Suppose that
	$a=0$ and $\widetilde a>0$. If
	$\mathcal F(0)=\mathcal F(\widetilde x)$, then
	$
	\pi z_{\widetilde a}
	\left(\frac{T(\widetilde a)}{2\pi}
	\widetilde\theta\right)
	=(A,0).
	$
	By the injectivity of $\pi$ on $\K$, we obtain
	\[
	z_{\widetilde a}
	\left(\frac{T(\widetilde a)}{2\pi}
	\widetilde\theta\right)
	=(A,0,B,0).
	\]
	Since $(A,0,B,0)$ is an equilibrium,
	uniqueness for \eqref{ODE} forces
	$z_{\widetilde a}$ to be constant,
	contradicting the construction of
	$\gamma_{\widetilde a}$.
	Thus \(a=0\) if and only if \(\widetilde a=0\). This proves that
	$\mathcal F$ is injective.

	By invariance of domain, $\mathcal F(\mathcal D)$ is open in $\mathbb R^2$, and $\mathcal F $ is a homeomorphism onto its image.
	Therefore
	\[
	\mathcal U
	:=
	\mathcal F(\mathcal D)
	=
	\{(A,0)\}
	\cup
	\bigcup_{0<a<a_*}\pi\gamma_a
	\]
	is an open neighborhood of $(A,0)$. This proves
	\textup{(i)}.

	\medskip\noindent\emph{Proof of \textup{(ii)}: Local continuation of periodic orbits.}
	Let $z_0(t)=(U(t),U'(t),V(t),V'(t))$ be a
	parametrization of $\gamma_0$ with least period $T_0$.
	The corresponding variational equations are
	\begin{equation}\label{hj1}
		h''+dh'-ch+pV^{p-1}j=0,
		\qquad
		j''-dj'-cj+qU^{q-1}h=0.
	\end{equation}
	
	Introducing the variables
	\begin{equation}\label{hj2}
		z_1=h,\qquad z_2=h'+\beta h,\qquad
		z_3=-j,\qquad z_4=-j'-\alpha j,
	\end{equation}
	we obtain
	\begin{equation}\label{positivevariation}
		\begin{aligned}
			z_1'&=-\beta z_1+z_2,&
			z_2'&=\alpha z_2+pV^{p-1}z_3,\\
			z_3'&=-\alpha z_3+z_4,&
			z_4'&=\beta z_4+qU^{q-1}z_1.
		\end{aligned}
	\end{equation}
	
	Writing $\mathbf z=(z_1,z_2,z_3,z_4)^{\mathsf T}$, we have
	\[
	\mathbf z'=B(t)\mathbf z,
	\qquad
	B(t)=
	\begin{pmatrix}
		-\beta & 1 & 0 & 0\\
		0 & \alpha & pV^{p-1} & 0\\
		0 & 0 & -\alpha & 1\\
		qU^{q-1} & 0 & 0 & \beta
	\end{pmatrix}.
	\]
	Since $U,V>0$, the matrix $B(t)$ is cooperative and irreducible.
	Thus,
	the evolution matrix \(\mathcal X(t,s)\) has strictly positive
	entries whenever \(t>s\). In particular,
	$
	M:=\mathcal X(T_0,0)
	$
	is a strictly positive matrix. By the Perron--Frobenius theorem, its
	spectral radius \(\mu=\rho(M)>0\) is an algebraically simple
	eigenvalue of \(M\), and \(M\) admits a strictly positive
	eigenvector associated with \(\mu\).

	Notice that  the tangent variation $(h,j)=(U',V')$ is a nonzero $T_0$-periodic
	solution of the variational system \eqref{hj1}. Under the change of variables in
	\eqref{hj2},
	we denote by $\mathbf z_{\mathrm{tan}}$ the corresponding $T_0$-periodic solution of
	\eqref{positivevariation}. Thus,
	we have
	\[
	M\mathbf z_{\mathrm{tan}}(0)
	=\mathbf z_{\mathrm{tan}}(T_0)=
	\mathbf z_{\mathrm{tan}}(0),
	\]
	which implies that $1$ is an eigenvalue of $M$, and therefore
	$\mu\ge1$.
	
	Suppose that $\mu=1$.
	Since $M$ is strictly positive, by the
	Perron--Frobenius theorem, it follows that the
	eigenspace corresponding to the eigenvalue $1$
	is one-dimensional and spanned by a strictly
	positive vector. Thus, $\mathbf z_{\mathrm{tan}}(0)$
	must be a nonzero scalar multiple of this vector.
	As $\mathcal X(t,0)$ has strictly positive entries
	for every $t>0$, the first component of
	$\mathbf z_{\mathrm{tan}}(t)$, namely $U'(t)$,
	has a fixed strict sign on $(0,T_0)$.
	This contradicts the periodicity of $U$, since
	\[
	\int_0^{T_0}U'(t)\,\dd t
	=
	U(T_0)-U(0)
	=
	0.
	\]
	Therefore, $\mu>1$.

	Finally, $M$ is similar to the monodromy matrix
	of the variational system in canonical Hamiltonian
	coordinates. Since the latter is symplectic,
	the eigenvalues of $M$ occur in reciprocal pairs,
	and $\det M=1$.
	As $\mu>1$ is an eigenvalue, so is $\mu^{-1}$.
	Moreover, the autonomy of the system ensures that
	$1$ is an eigenvalue. Since $M$ is four-dimensional
	and its determinant is one, the remaining eigenvalue
	must also equal $1$. Consequently,
	\[
	\operatorname{spec}(M)
	=
	\{\mu,\mu^{-1},1,1\},
	\qquad \mu>1,
	\]
	where eigenvalues are counted with their algebraic
	multiplicities.
	
	We now verify that the energy surface
	$
	\Sigma_{E_0}:=\{\cH=E_0\}
	$
	is regular near $\gamma_0$. In fact,
	direct differentiation of
	\eqref{Hamiltonian} shows that
	every critical point of $\cH$ is an
	equilibrium of \eqref{ODE}. 
	By uniqueness for the initial-value
	problem, a trajectory that passes through an equilibrium must be
	constant. Since $\gamma_0$ is a nonconstant periodic orbit, it follows
	that
	\[
	\dd\cH_z\ne0
	\qquad\text{for every }z\in\gamma_0.
	\]
	Since $\gamma_0$ is compact and $\dd\cH$ is continuous, there exists
	an open neighborhood $\mathcal N$ of $\gamma_0$ such that
	\[
	\dd\cH_z\ne0
	\qquad\text{for every }z\in\mathcal N.
	\]
	It follows from the implicit function theorem that
	$\Sigma_{E_0}\cap\mathcal N$ is a smooth three-dimensional
	hypersurface.
	
	Since the Floquet multiplier $1$
	has algebraic multiplicity exactly two, by the Poincar\'e continuation
	result \cite[Chapter~4, Section~4.1, Proposition~2]{HZ},
	there exist $\delta>0$ and a smooth family of
	periodic orbits $\gamma_E$, with  periods $T(E)$,
	such that
	\[
	\cH|_{\gamma_E}=E,\qquad
	\gamma_{E_0}=\gamma_0,\qquad
	T(E_0)=T_0,
	\qquad |E-E_0|<\delta.
	\]
	Choose a local three-dimensional section $\mathcal D$ through $z_0(0)$ transverse
	to the flow. For each $E$, we choose the intersection
	of $\gamma_E$ with $\mathcal D$ near $z_0(0)$ as the initial point,
	and denote the corresponding solution of \eqref{ODE} by $z_E(t)$.
	Then $z_E(t)$ depends smoothly on $(E,t)$ and
	\[
	z_E(\mathbb R)=\gamma_E,\qquad
	z_E(t+T(E))=z_E(t).
	\]
	Since $z_{E_0}(0)=z_0(0)$, by the uniqueness for the initial-value problem, we have
	\[
	z_{E_0}(t)=z_0(t)
	\qquad\text{for all }t\in\mathbb R.
	\]

	Write
	$
	z_E=(U_E,U_E',V_E,V_E').
	$
	Since $U$ and $V$ are positive and $T_0$-periodic, then
	$
	m_0:=\min_{t\in[0,T_0]}\min\{U(t),V(t)\}>0.
	$
	Noting that $z_E$ and $T(E)$ depend smoothly on $E$, we obtain 
	\[
	\sup_{\theta\in[0,1]}
	\left(
	|U_E(\theta T(E))-U(\theta T_0)|
	+
	|V_E(\theta T(E))-V(\theta T_0)|
	\right)\longrightarrow0
	\qquad\text{as }E\to E_0.
	\]
	After decreasing $\delta$ if necessary, it follows that
	\[
	U_E(t)\ge\frac{m_0}{2},
	\qquad
	V_E(t)\ge\frac{m_0}{2}
	\qquad
	\text{for all }t\in[0,T(E)]
	\text{ and }|E-E_0|<\delta.
	\]
	By periodicity, these inequalities hold for every $t\in\mathbb R$.
	Consequently, every
	$\gamma_E$ is a positive periodic orbit.

	We next show that, for $\delta$ sufficiently small, $T(E)$ is the least period
	of $z_E$. Otherwise, there exist $E_j\to E_0$ such that the least
	period $S_j$ of $z_{E_j}$ satisfies
	$
	0<S_j\le \frac{T(E_j)}2.
	$
	Passing to a subsequence, we have
	\[
	S_j\to S\in[0,T_0/2], \quad \text{as} \ j\rightarrow +\infty.
	\]
	If $S=0$, by the continuity of $\dot z_{E_j}(t)$ with respect to $(E,t)$, we yield 
	\[
	0=
	\frac{z_{E_j}(S_j)-z_{E_j}(0)}{S_j}
	=
	\int_0^1\dot z_{E_j}(\theta S_j)\,\dd\theta
	\rightarrow \dot z_0(0), \quad \text{as} \ j\rightarrow +\infty,
	\]
	contradicting $\dot z_0(0)\ne0$. Thus $S>0$. Since  $z_{E_j}(S_j)=z_{E_j}(0)$, letting $j\rightarrow +\infty$, we have 
	\[
	z_0(S)=z_0(0),
	\qquad 0<S\le T_0/2,
	\]
	contradicting the minimality of $T_0$.
	Thus $T(E)$ is the least
	period of $z_E$ for every $|E-E_0|<\delta$.
	
	Now define
	\[
	F:(E_0-\delta,E_0+\delta)\times(\mathbb R/\mathbb Z)
	\longrightarrow\mathbb R^2
	\]
	by
	\[
	F(E,\theta):=\pi z_E\bigl(\theta T(E)\bigr).
	\]
	The map $F$ is continuous. It is also injective. Indeed, suppose that
	$
	F(E,\theta)=F(\widetilde E,\widetilde\theta).
	$
	By the injectivity of $\pi$ on $\K$, established in
	Lemma~\ref{lem:projection}, we have 
	$
	z_E\bigl(\theta T(E)\bigr)
	=
	z_{\widetilde E}
	\bigl(\widetilde\theta T(\widetilde E)\bigr).
	$
	Thus,  
	$
	E
	=\widetilde E,
	$
	and  $$
	z_E\bigl(\theta T(E)\bigr)
	=
	z_{ E}
	\bigl(\widetilde\theta T( E)\bigr).
	$$
	Since $T(E)$ is its least period, it follows that
	$
	\theta=\widetilde\theta$ in 
	$\mathbb R/\mathbb Z,
	$
	and hence $F$ is injective.
	
	Similarly, by invariance of domain, $F$ is a homeomorphism
	onto an open subset of $\mathbb R^2$. Consequently,
	\[
	\mathcal A
	:=
	\bigcup_{|E-E_0|<\delta}\pi\gamma_E
	\]
	is an open annulus containing $\pi\gamma_0$.
	This completes the proof of \textup{(ii)}.	
\end{proof}

\subsection{A criterion for periodicity}
\label{sec:periodicitycriterion}

We next establish a periodicity criterion using the planar
reduction and the local neighborhoods constructed in
Proposition~\ref{localperiodic}. When a limit set avoids
the origin, the Poincar\'e--Bendixson theorem applies
to the projected flow, and the injectivity of $\pi$
allows us to recover the corresponding periodic orbit
of the original system. This criterion plays a key role
in the subsequent global continuation argument.

\Needspace{7\baselineskip}
\begin{proposition}
	\label{periodicitycriterion}
	Let $z:\R\to\K\setminus\{0\}$ be a trajectory. If either its
	$\alpha$-limit set or its $\omega$-limit set does not contain the
	origin, then $z$ is constant or periodic.
\end{proposition}

\begin{proof}
	It suffices to consider the $\omega$-limit set, since the
	$\alpha$-limit case follows by reversing time.  Set
	\[
	\Omega:=\omega(z)
	=\bigcap_{s>0}\overline{\{z(t):t\ge s\}},
	\qquad 0\notin\Omega.
	\]
	Since $z(\mathbb R)\subset\mathcal K\setminus\{0\}$ and $\mathcal K$
	is compact by Lemma~\ref{Kcompact}, the positive orbit of $z$ is
	relatively compact in $\mathcal K$. Moreover, $\mathcal K$ is invariant
	under the flow of \eqref{ODE} by its definition. Hence, $\Omega=\omega(z)$ is nonempty, compact, and invariant.
	Moreover,
	\[
	\operatorname{dist}\bigl(z(t),\Omega\bigr)\rightarrow0
	\qquad\text{as }t\to+\infty.
	\]
	Since $0\notin\Omega$ and $\Omega$ is compact, then 
	$
	d_0:=\operatorname{dist}(0,\Omega)>0.
	$
	It follows that  for some $t_0>0$,
	\[
	\operatorname{dist}\bigl(z(t),\Omega\bigr)<\frac{d_0}{2}
	\qquad (t\ge t_0),
	\]
	from which we have 
	\[
	|z(t)|
	\ge \operatorname{dist}(0,\Omega)
	-\operatorname{dist}(z(t),\Omega)
	>\frac{d_0}{2}
	\qquad (t\ge t_0).
	\]
	Thus
	$
	C:=\overline{\{z(t):t\ge t_0\}}
	$
	is a compact subset of $\mathcal K\setminus\{0\}$.
	
	Since 
	$\pi:\mathcal K\to\pi\mathcal K$ is a homeomorphism by
	Lemma~\ref{lem:projection}, 
	$
	K_0:=\pi(C)
	$
	is a compact subset of $\pi\mathcal K\setminus\{0\}$.
	Applying Lemma~\ref{lem:planar} to $K_0$, there exists an open
	neighborhood $\mathcal O$ of $\pi(C)$ and a globally Lipschitz vector
	field
	$
	\widetilde{\mathcal F}:\mathbb R^2\to\mathbb R^2
	$
	such that
	$
	\widetilde{\mathcal F}=\mathcal F$ on
	$\mathcal O\cap\pi\mathcal K.
	$
	Let
	\[
	x(t):=\pi z(t).
	\]
	Since $x(t)\in\pi(C)\subset\mathcal O$ for $t\ge t_0$, the projected
	orbit satisfies
	\[
	x'(t)=\widetilde{\mathcal F}(x(t))
	\qquad(t\ge t_0).
	\]
	Thus the positive tail of $x$ is a trajectory of the globally
	Lipschitz planar system
	$
	\dot x=\widetilde{\mathcal F}(x),x\in\R^2.
	$
	Moreover, since $\pi|_C:C\to\pi(C)$ is a homeomorphism, it preserves
	the corresponding limit sets. Hence
	\[
	\omega(x)
	=\omega(\pi z)
	=\pi\omega(z)
	=\pi\Omega.
	\]
	
	Suppose first that
	$
	z_*:=(A,0,B,0)\in\Omega.
	$
	By the definition of $\Omega=\omega(z)$, there exists a sequence
	$t_j\to+\infty$ such that
	$
	z(t_j)\rightarrow z_*.
	$
	Hence, by the continuity of $\pi$, we have 
	\[
	x(t_j)=\pi z(t_j)\rightarrow\pi z_*=(A,0).
	\]
	By Proposition~\ref{localperiodic}\textup{(i)}, the set
	$
	\mathcal U
	:=
	\{(A,0)\}
	\cup
	\bigcup_{0<a<a_*}\pi\gamma_a
	$
	is an open neighborhood of $(A,0)$. Therefore,
	$
	x(t_j)\in\mathcal U
	$
	for all sufficiently large $j$.
	
	Suppose next that
	$
	z_*:=(A,0,B,0)\notin\Omega.
	$
	Since $\pi$ is injective, it follows that
	$
	(A,0)\notin\pi\Omega.
	$
	By Lemma~\ref{lem:planar}, $(A,0)$ is the unique equilibrium of the
	projected vector field in $\pi\mathcal K\setminus\{0\}$. Hence the
	compact $\omega$-limit set
	$
	\omega(x)=\pi\Omega
	$
	contains no equilibrium. Applying the Poincar\'e--Bendixson theorem
	\cite[Chapter~6, Theorem~15]{HO}, we obtain 
	$
	\pi\Omega=\Gamma,
	$
	where $\Gamma$ is a periodic orbit of the extended planar system
	$
	\dot x=\widetilde{\mathcal F}(x),x\in\R^2.
	$

	We next  show that $\Omega$ is a positive
	periodic orbit of \eqref{ODE}. In fact,
	choose $\zeta\in\Omega$, and let $Z(t)$ be the solution of
	\eqref{ODE} with $Z(0)=\zeta$. By the invariance of $\Omega$,
	$
	Z(\mathbb R)\subset\Omega.
	$
	Since $X(t):=\pi Z(t)$ is a trajectory on the periodic orbit
	$\Gamma=\pi\Omega$, there exists $T>0$ such that
	\[
	\pi Z(T)=X(T)=X(0)=\pi Z(0).
	\]
	By 
	Lemma~\ref{lem:projection}, we have  $Z(T)=Z(0)$, and hence $Z$ is periodic.
	Moreover,
	\[
	\pi\bigl(Z(\mathbb R)\bigr)=\Gamma=\pi\Omega,
	\]
	and by the injectivity of $\pi$ again, we yield
	$Z(\mathbb R)=\Omega$. Since $0\notin\Omega$, $\Omega$ is a positive
	periodic orbit of \eqref{ODE}.
	
	Applying Proposition~\ref{localperiodic}\textup{(ii)} to the periodic
	orbit $\Omega$, there exists an open annulus $\mathcal A_\Omega$
	containing
	$
	\Gamma=\pi\Omega.
	$
	Since $\Gamma$ is compact, there exists $\delta_\Omega>0$ such that
	\[
	\{y\in\mathbb R^2:
	\operatorname{dist}(y,\Gamma)<\delta_\Omega\}
	\subset\mathcal A_\Omega.
	\]
	Moreover, we have 
	\[
	\operatorname{dist}(x(t),\Gamma)
	=
	\operatorname{dist}\bigl(\pi z(t),\pi\Omega\bigr)
	\le
	\operatorname{dist}\bigl(z(t),\Omega\bigr)
	\rightarrow0,\quad \text{as} \ t\rightarrow \infty.
	\]
	Consequently,
	$
	x(t)\in\mathcal A_\Omega
	$
	for all sufficiently large $t$.
	
	In either case, there exist $s\in\mathbb R$ and
	a constant or positive periodic solution
	$\widetilde z$ of \eqref{ODE} such that
	\[
	\pi z(s)=\pi\widetilde z(0),
	\]
	which implies
	\[
	z(s)=\widetilde z(0).
	\]
	Since $w(t):= z(s+t)$ and $\widetilde z(t)$ are the solutions to \eqref{ODE} with $w(0)=\widetilde z(0)$,  we have
	\[
	z(s+t)=\widetilde z(t)
	\qquad\text{for all }t\in\mathbb R.
	\]
	Consequently, $z$ is either constant or periodic.
\end{proof}

\subsection{Global continuation and parametrization}
\label{sec:globalperiodic}

We next extend the local families of periodic orbits globally,
parametrizing them by the minimum of the first component.
Combining uniqueness, compactness, and the preceding periodicity
criterion, we show that the parameter covers the entire interval
$(0,A)$. We further establish the energy range, characterize
the critical points of both components, and determine the
limiting behavior of the periods.

\Needspace{7\baselineskip}

\begin{proposition}[Global parametrization of positive periodic orbits]
	\label{periodicfamily}
	For every $\eps\in(0,A)$, there exists a unique positive nonconstant
	periodic orbit $\gamma_\eps\subset\mathcal K$ such that
	\[
	\min_{t\in\mathbb R}U_\eps(t)=\eps.
	\]
	Then the energy map
	\[\mathcal{E}
	:(0,A)\longrightarrow(E_*,0),
	\qquad
	\mathcal{E}(\eps):=\mathcal H(z),
	\qquad z\in\gamma_\eps.
	\]
	is a continuous strictly decreasing bijection.
	The corresponding components $U_\eps$ and $V_\eps$ have the same
	least period $T_\eps$ satisfying \eqref{periodlimits}, and have exactly one
	nondegenerate minimum and one nondegenerate maximum per least period.
\end{proposition}

\begin{proof}
	\medskip\noindent\emph{Step 1: Extremal bounds and uniqueness.}
	Let $\gamma\subset\K$ be a positive nonconstant periodic orbit, and
	let
	\[
	z(t)=(U(t),U'(t),V(t),V'(t)),
	\qquad z(\R)=\gamma.
	\]
	Set
	\[
\begin{aligned}
 m&:=\min_{\R}U=U(t_m), &\quad \ell&:=\min_{\R}V=V(t_\ell),\\
 M&:=\max_{\R}U=U(t_M), &\quad L&:=\max_{\R}V=V(t_L).
\end{aligned}
\]
It follows by \eqref{ODE} that 
	\begin{align*}
		0
		&\le U''(t_m)
		=cm-V(t_m)^p
		\le cm-\ell^p,\\
		0
		&\le V''(t_\ell)
		=c\ell-U(t_\ell)^q
		\le c\ell-m^q,\\
		0
		&\ge U''(t_M)
		=cM-V(t_M)^p
		\ge cM-L^p,\\
		0
		&\ge V''(t_L)
		=cL-U(t_L)^q
		\ge cL-M^q.
	\end{align*}
	Consequently, we have
	\begin{align*}
		m^{pq-1}
		&\le c^{p+1}=A^{pq-1},
		&
		\ell^{pq-1}
		&\le c^{q+1}=B^{pq-1},\\
		M^{pq-1}
		&\ge c^{p+1}=A^{pq-1},
		&
		L^{pq-1}
		&\ge c^{q+1}=B^{pq-1}.
	\end{align*}
	Moreover, if $m=A$ or $M=A$, then $(U,U')=(A,0)$ at $t_{m}$ or $t_{M}$.
	Similarly, if $\ell=B$ or $L=B$, then $(V,V')=(B,0)$ at $t_{\ell}$ or $t_{L}$.
	Thus, we obtain from Lemma~\ref{lem:uniqueness} that 
	\begin{equation}\label{strict-minmax}
		0<m<A<M,\qquad 0<\ell<B<L.
	\end{equation}
	
	Let $z_1,z_2$ be positive nonconstant periodic solutions such that
	$\min_{\R}U_1=\min_{\R}U_2=m$. 
	We may assume that
	\[
	U_1(0)=U_2(0)=m,
	\qquad
	U_1'(0)=U_2'(0)=0.
	\]
	Using Lemma~\ref{lem:uniqueness}, we have 
	\[
	z_1(t)=z_2(t)
	\qquad(t\in\R).
	\]
	Similarly, if
	\[
	U(t_1)=U(t_2)=m,
	\qquad
	U'(t_1)=U'(t_2)=0,
	\]
	then using Lemma~\ref{lem:uniqueness} again, it follows that 
	\[
	z(t+t_1)=z(t+t_2)
	\qquad(t\in\R).
	\]
	Consequently, a positive nonconstant periodic orbit is uniquely
	determined by $\min_{\R}U$, and $\min_{\R}U$ is
	attained at a unique phase modulo the least period.
	
	\medskip\noindent\emph{Step 2: Global continuation by the minimum.}
	Define
	\[
	S:=
	\left\{
	\min_{t\in\R}U(t):
	z(t)
	\text{ is a positive nonconstant periodic solution with }
	z(\R)\subset\K
	\right\} \subset(0,A).
	\]
	Fix $\eps_0\in S$. By Step~1, there is a unique periodic orbit
	$\gamma_{\eps_0}\subset\K$ such that
	$
	\min_{t\in\R}U_{\eps_0}(t)=\eps_0.
	$
	Let
	$
	E_0:=\cH|_{\gamma_{\eps_0}},
	$
	and let $\{\gamma_E:|E-E_0|<\delta\}$ be the local family
	from Proposition~\ref{localperiodic}\textup{(ii)}.
	Choose periodic representatives $z_E$ and periods $T(E)$
	depending smoothly on $E$, as in its proof, and let $U_E$
	denote the first component of $z_E$. Define
	\begin{align*}
	\mu(E)
	:= & ~ \min_{t\in\R}U_E(t)
	= \min_{\theta\in[0,1]}U_E(\theta T(E)).
	\end{align*}
	For $E,\widetilde E\in(E_0-\delta,E_0+\delta)$, we have
	\begin{align*}
	|\mu(E)-\mu(\widetilde E)|
	\leq & ~
	\sup_{\theta\in[0,1]}
	\left|
	U_E(\theta T(E))
	-
	U_{\widetilde E}(\theta T(\widetilde E))
	\right|.
	\end{align*}
	The smooth dependence of the periodic profiles and periods
	therefore implies that $\mu$ is continuous.

	Moreover, $\mu$ is injective. Indeed, if
	$\mu(E_1)=\mu(E_2)$, by Step~1, we have 
	$\gamma_{E_1}=\gamma_{E_2}$, and evaluating $\cH$ on this
	common orbit yields $E_1=E_2$.
	Thus $\mu$ is strictly monotone, and
	\begin{align*}
	\mu\bigl((E_0-\delta,E_0+\delta)\bigr)
	\end{align*}
	is an open interval contained in $S$ and containing
	$\mu(E_0)=\eps_0$. Since $\eps_0\in S$ was arbitrary,
	$S$ is open.
	Finally, using Proposition~\ref{localperiodic}\textup{(i)}, there exists a 
	$\eta>0$ such that
	$
	(A-\eta,A)\subset S.
	$

	Let $(a_0,A)$ be the  connected component of $S$ containing $(A-\eta,A)$, where \(0\le a_0<A\).
	If $a_0>0$, choose $\left\{\eps_j\right\}\subset(a_0,A)$ such that 
	\[
	\eps_j\rightarrow a_0 \quad (\text{as}~ j\rightarrow\infty)
	,\quad
	U_j(0)=\eps_j,\qquad U_j'(0)=0.
	\]
	By Lemma~\ref{Kcompact}, up to a subsequence, there exists $ z_0\in\mathcal K\setminus\{0\}$ such that
	\[
	\begin{aligned}
		z_j(0)\rightarrow z_0,\quad
		z_j\rightarrow z=(U,U',V,V')
		\quad\text{in }C^1_{\mathrm{loc}}(\R), \quad \text{as}~ j\rightarrow\infty.
	\end{aligned}
	\]
	Letting $j\rightarrow\infty$, we have 
	\[
	U(0)=a_0,
	\qquad
	U(t)\ge a_0\quad(t\in\R).
	\]
	Consequently,
	$
	0\notin\alpha(z)\cup\omega(z). 
	$
	It follows by Proposition~\ref{periodicitycriterion} that $z$ is constant or
	periodic. Since $U(0)=a_0<A$, together with \eqref{strict-minmax},  $z$ is periodic and
	$
	\min_{\R}U=a_0.
	$
	Thus $a_0\in S$. Since $S$ is open, there exists $\delta_0>0$ such
	that
	$
	(a_0-\delta_0,a_0+\delta_0)\subset S,
	$
	contradicting the definition of $a_0$.
	Hence
	$
	S=(0,A).
	$
	
	Therefore, for every $\eps\in(0,A)$, there exists a positive
	nonconstant periodic orbit $\gamma_\eps\subset\K$ satisfying
	$
	\min_{t\in\R}U_\eps(t)=\eps.
	$
	By Step~1, this orbit is unique. Hence,
	$
	\{\gamma_\eps:0<\eps<A\}
	$
	is precisely the family of all positive nonconstant periodic orbits
	in $\K$.

	\medskip\noindent\emph{Step 3: The energy parametrization.}
	For each $\eps\in(0,A)$, let
	$
	E_\eps:=\cH|_{\gamma_\eps}
	$
	denote the conserved energy of the orbit
	$\gamma_\eps$.
	For fixed  $\eps_0\in(0,A)$, by Step~2, the map
	\[
	\mu(E):=\min_{t\in\R}U_E(t),
	\qquad |E-E_0|<\delta,
	\]
	is continuous and strictly monotone, with $\mu(E_0)=\eps_0$.
	For every $E\in(E_0-\delta,E_0+\delta)$, the periodic orbits
	$\gamma_E$ and $\gamma_{\mu(E)}$ satisfy 
	$\min_{\mathbb R}U_E
	=
	\mu(E)
	=
	\min_{\mathbb R}U_{\mu(E)}.$
	Hence by Step~1, we have  
	$
	\gamma_E=\gamma_{\mu(E)}.
	$
	Since $\gamma_E$ has energy $E$, it follows that
	$
	E_{\mu(E)}=E.
	$
	Therefore, on the open interval
	$
	\mu\bigl((E_0-\delta,E_0+\delta)\bigr),
	$
	the map $\eps\mapsto E_\eps$ is the inverse of $\mu$. In particular,
	this map is continuous and injective on
	the neighborhood $I_{\eps_0}$ of $\eps_0$.Thus, 
	$
	\eps\mapsto E_\eps
	$
	is either strictly increasing or strictly decreasing on $(0,A)$.

	For each $\eps\in(0,A)$, set
	$
	z_\eps=(U_\eps,U_\eps',V_\eps,V_\eps')
	$ with $U_\eps$ attaining its minimum at $t=0$. Then
	$
	z_\eps(0)
	=
	\bigl(\eps,0,V_\eps(0),V_\eps'(0)\bigr).
	$
	It follows by Lemma~\ref{radialbounds} that 
	\[
	0<V_\eps(0)
	\le \kappa^{-1/p}\eps^{1/p},
	\qquad
	|V_\eps'(0)|
	\le \max\{\alpha,\beta\}V_\eps(0).
	\]
	Consequently, we have 
	\begin{equation}\label{energy-limit-zero}
		z_\eps(0)\rightarrow0,
		\qquad
		E_\eps=\cH(z_\eps(0))
		\rightarrow\cH(0)=0
		\qquad\text{as} ~\eps\rightarrow0.
	\end{equation}
	
	Recall from the proof of
	Proposition~\ref{localperiodic}\textup{(i)} that
	\begin{align*}
	m(a)
	:= & ~ \min_{t\in\R}U_a(t)
	= A-a+O(a^2),
	\qquad
	m'(a)<0,
	\end{align*}
	and $m((0,a_*))=(A-\eta,A)$.
	For $\eps\in(A-\eta,A)$, define
	$a(\eps):=m^{-1}(\eps)$.
	Then $a(\eps)\to0$ as $\eps\to A$.
	Since $\gamma_\eps$ and $\gamma_{a(\eps)}$ have the same
	minimum of their first components, Step~1 gives
	\begin{equation}\label{eq:small-amplitude-identification}
	\gamma_\eps=\gamma_{a(\eps)}
	\qquad (A-\eta<\eps<A).
	\end{equation}
	The uniform convergence of the small-amplitude orbits to
	$(A,0,B,0)$ therefore implies, by continuity of $\cH$, that
	\begin{equation}\label{energy-limit-equilibrium}
	z_\eps(0)\longrightarrow(A,0,B,0),
	\qquad
	E_\eps\longrightarrow E_*
	\qquad\text{as }\eps\to A.
	\end{equation}
	
	Since $E_*<0$, we deduce from \eqref{energy-limit-zero} and
	\eqref{energy-limit-equilibrium} that the map
	$\eps\mapsto E_\eps$ is strictly decreasing, and
	$
	\{E_\eps:0<\eps<A\}=(E_*,0).
	$
	Therefore,
	\[
	(0,A)\ni\eps\longmapsto E_\eps\in(E_*,0)
	\]
	is a continuous strictly decreasing bijection.
	
	\medskip\noindent\emph{Step 4: Periods and critical points.}
	If $U_\eps$ is constant, by \eqref{ODE}, we have  
	$V_\eps^p\equiv cU_\eps$, and hence $V_\eps$ is constant. Similarly,
	if $V_\eps$ is constant, we also obtain that $U_\eps$
	is constant.
	Thus both $U_\eps$ and $V_\eps$ are nonconstant. For a nonconstant periodic
	function $W$, set
	\[
	\mathcal P(W):=\{P>0:W(t+P)=W(t)\ \text{for all }t\in\R\}.
	\]
	Let $P\in\mathcal P(U_\eps)$. It follows from \eqref{ODE} that 
	\[
	V_\eps^p=cU_\eps-U_\eps''-dU_\eps'
	\]
	is $P$-periodic, and therefore $P\in\mathcal P(V_\eps)$. Interchanging $U_\eps$ and $V_\eps$, 
	we lead
	\[
	\mathcal P(U_\eps)=\mathcal P(V_\eps),
	\qquad
	T_\eps:=\min\mathcal P(U_\eps)
	=\min\mathcal P(V_\eps).
	\]
	
	We next count the critical points in one least period. Define
	\[
	\begin{aligned}
		N_U(\eps)
		&:=
		\#\{t\in\R/(T_\eps\mathbb Z):U_\eps'(t)=0\},\\
		N_V(\eps)
		&:=
		\#\{t\in\R/(T_\eps\mathbb Z):V_\eps'(t)=0\}.
	\end{aligned}
	\]
	Applying  Lemma~\ref{lem:uniqueness}, we have 
	\[
	U_\eps''(t)\ne0\quad\text{whenever }U_\eps'(t)=0,
	\qquad
	V_\eps''(t)\ne0\quad\text{whenever }V_\eps'(t)=0.
	\]

	By Step~3, the local families in
	Proposition~\ref{localperiodic}\textup{(ii)} may be parametrized by
	$\eps$. 
	Define
	\[
	\widehat U_\eps(\theta):=U_\eps(T_\eps\theta),
	\qquad
	\widehat V_\eps(\theta):=V_\eps(T_\eps\theta),
	\qquad
	\theta\in\mathbb R/\mathbb Z.
	\]
	By Proposition~\ref{localperiodic}\textup{(ii)} and Step~3,
	$\widehat U_\eps$ and $\widehat V_\eps$ depend continuously
	on $\eps$ in $C^2(\mathbb R/\mathbb Z)$.
	Since all their critical points are nondegenerate by
	Lemma~\ref{lem:uniqueness},
	$N_U$ and $N_V$ are locally constant on $(0,A)$.
	For $\eps$ sufficiently close to $A$,
	Eq.~\eqref{eq:small-amplitude-identification} and
	Proposition~\ref{localperiodic}\textup{(i)} give
	$N_U(\eps)=N_V(\eps)=2$.
	Since $(0,A)$ is connected, it follows that
	\[
	N_U(\eps)=N_V(\eps)=2
	\qquad(0<\eps<A).
	\]
	Thus each component has one nondegenerate minimum and one
	nondegenerate maximum in $\R/(T_\eps\mathbb Z)$.
	
	It remains to prove the period limits. by Lemma~\ref{radialbounds}
	\ref{derivativebounds}, we obtain  
	\[
	\left|(\log U_\eps)'(t)\right|
	\le M_0:=\max\{\alpha,\beta\}.
	\]
	If $t_-$ and $t_+$ are the minimum and maximum points of $U_\eps$
	in one least period, together with \eqref{strict-minmax}, we  yield
	\[
	\log\frac{A}{\eps}
	<
	\log\frac{U_\eps(t_+)}{U_\eps(t_-)}
	=
	\left|\int_{t_-}^{t_+}\frac{U_\eps'(t)}{U_\eps(t)}\,\dd t\right|
	\le M_0T_\eps.
	\]
	Hence
	\[
	T_\eps\rightarrow+\infty
	\qquad\text{as} ~\eps\rightarrow0.
	\]
	Finally, Eq.~\eqref{eq:small-amplitude-identification}
	identifies the two periodic orbits and hence their least periods.
	Since $a(\eps)\to0$ as $\eps\to A$,
	Proposition~\ref{localperiodic}\textup{(i)} yields
	\begin{align*}
	T_\eps
	= & ~ T(a(\eps))
	\longrightarrow \frac{2\pi}{\omega}
	\qquad\text{as }\eps\to A.
	\end{align*}
\end{proof}

\subsection{Asymptotic behavior of nonperiodic trajectories}
\label{sec:nonperiodic}
We finally consider positive global trajectories $z(t)$ that are
neither constant nor periodic. By the preceding periodicity
criterion, both limit sets must contain the origin.
Combining monotonicity below the equilibrium level with
the component estimates, we prove that $z(t)\to0$ as
$t\to\pm\infty$ and determine the corresponding energy.

\begin{proposition}
	\label{remaining}
	Let $z=(U,U',V,V')$ be a positive global trajectory that is
	neither constant nor periodic. Then
	\[
	\lim_{t\to-\infty}z(t)
	=
	\lim_{t\to+\infty}z(t)
	=
	0,
	\qquad
	\cH(z(t))=0
	\quad(t\in\R).
	\]
\end{proposition}

\begin{proof}

	We first claim that
	\begin{equation}\label{noturning}
		U'(t)\ne0
		\qquad\text{whenever }U(t)\le A.
	\end{equation}
	Suppose that $U'(t_0)=0$. Since $U>0$, there are two cases.
	If
	$
	\eps_0:=U(t_0)\in(0,A),
	$
	let
	$
	z_{\eps_0}
	=
	(U_{\eps_0},U_{\eps_0}',V_{\eps_0},V_{\eps_0}')
	$
	parametrize the periodic orbit $\gamma_{\eps_0}$, with its minimum
	placed at $t=0$. Then
	\[
	U_{\eps_0}(0)=\eps_0=U(t_0),
	\qquad
	U_{\eps_0}'(0)=0=U'(t_0).
	\]
	By Lemma~\ref{lem:uniqueness}, $z$ is periodic, a contradiction.
	If $U(t_0)=A$, then $z$ and the positive equilibrium
	$
	z_*
	$
	satisfy
	$$
	(U(t_0),U'(t_0))=(A,0).
	$$
	Together with Lemma~\ref{lem:uniqueness} again, then 
	$z$ is constant, a contradiction. This proves
	\eqref{noturning}.

	Set
	\[
	\mathcal O:=\{t\in\mathbb R:U(t)<A\}.
	\]
	We obtain from  \eqref{noturning} that  $U$ is strictly monotone on every connected
	component of $\mathcal O$. Suppose that $(a,b)$ is a bounded
	connected component of $\mathcal O$. By the continuity of $U(t)$, we have 
	\[
	U(a)=U(b)=A,
	\]
	contrary to the strict monotonicity of $U$. Hence, $\mathcal O$ has no bounded connected component.

	Applying  Proposition~\ref{periodicitycriterion},
	we know that 
	$
	0\in\alpha(z)\cap\omega(z).
	$
	Then there exist sequences
	$$
	t_j^-\rightarrow-\infty,
	\quad
	t_j^+\rightarrow+\infty
	$$
	such that
	$$
	z(t_j^-)\rightarrow0,\quad 
	z(t_j^+)\rightarrow0$$
	as $j\rightarrow +\infty.$ Thus,
	for all sufficiently large $j$, we have 
	$t_j^-,t_j^+\in\mathcal O$, and hence
	$\mathcal O$ has a component $I_-$ unbounded to the left and a component $I_+$
	unbounded to the right,  
	these components are distinct. Otherwise, $I_-=I_+=\mathbb R$, and \eqref{noturning} implies that
	$U$ is strictly monotone on $\mathbb R$. Since $U(t_j^\pm)\to0$ as $j\to\infty$, the monotonicity of $U$ implies that $U\equiv0$, contradicting $U>0$. Therefore,
	for some $a\le b$,
	$
	I_-=(-\infty,a)$,
	and 
	$ I_+=(b,\infty).
	$
	Since $U$ is monotone on $I_-$ and $I_+$, we have 
	\[
	\lim_{t\to-\infty}U(t)
	=
	\lim_{t\to+\infty}U(t)
	=
	0.
	\]
	
	By Lemma~\ref{radialbounds}\ref{derivativebounds}  and \ref{crossbounds}, we have 
	\[
	\lim_{t\to\pm\infty}V(t)=0.
	\]
	and 
	\[
	\lim_{t\to\pm\infty}U'(t)
	=
	\lim_{t\to\pm\infty}V'(t)
	=
	0.
	\]
	Consequently,
	\[
	\lim_{t\to\pm\infty}z(t)=0.
	\]
	Finally, since $\cH$ is continuous and conserved along trajectories, we obtain 
	\[
	\cH(z(t))
	=
	\lim_{s\to+\infty}\cH(z(s))
	=
	\cH(0)
	=
	0
	\qquad(t\in\mathbb R).
	\]
\end{proof}

\section{Removability and completion of the classification}\label{sec:completion}
In this section, we establish a removability criterion for positive
radial solutions and complete the proof of
Theorem~\ref{classification}. Furthermore, we characterize the ground-state
homoclinic orbit and derive annular critical energy estimates
for removable and nonremovable solutions.

\subsection{Removability and the singular classification}

We first derive decay estimates that ensure
removability at the origin.

\begin{lemma}\label{removable}
Assume \eqref{critical} holds. Let $(U,V)$ be a positive global solution of \eqref{ODE} satisfying
	$(U(t),V(t))\to(0,0)$ as $t\to-\infty$. Then
	\begin{equation}\label{sharpdecay}
		U(t)=O(e^{\alpha t}),\qquad
		V(t)=O(e^{\beta t})
		\quad\text{as }t\to-\infty.
	\end{equation}
Moreover, the corresponding radial pair $(u,v)$
defined by \eqref{EF} extends to a positive
$C^2$ solution of \eqref{system} in $\mathbb R^n$.
\end{lemma}

\begin{proof}
	Interchanging the two components if necessary, we assume that
	$
	p\le q.
	$
	By the critical relation,
	$
	q\ge \frac{n+2}{n-2}>1,
	$
	whereas \(p\) is allowed to be sublinear. Since $p\beta=\alpha+2$ and $q\alpha=\beta+2$, we have
$
\alpha/(pq)<\beta/q<\alpha.
$
Fix $\sigma$ such that
\[
\frac{\alpha}{pq}<\sigma<\frac{\beta}{q},
\] and define
	\[
	C_\sigma=
	\frac{1}{(\alpha-\sigma)(\beta+\sigma)}
	\left[
	\frac1{(\beta-q\sigma)(\alpha+q\sigma)}
	\right]^{1/q}.
	\]
	Since $pq>1$, we may choose $\delta>0$ sufficiently
	small such that
	\begin{equation*}
	pq\delta^{p-1/q}
	\max\left\{c^{-1-1/q},C_\sigma\right\}\le\frac{1}{2}.
	\end{equation*}
	Since $(U,V)\to(0,0)$ as $t\to-\infty$, there exists
	$T<0$ such that
	\[
	0<U,V\le\delta\qquad\text{on }(-\infty,T].
	\]
	Since the system \eqref{ODE} is autonomous, after a translation,
	we may assume that $T=0$.

	The Green kernels corresponding to
	\[
	L_1=-D^2-dD+c,\qquad
	L_2=-D^2+dD+c
	\]
	on $(-\infty,0]$, with zero boundary value at $s=0$
	and boundedness as $s\to-\infty$, are given by
	\begin{equation*}
		\begin{aligned}
			G_1(s,t)
			&=\frac1{\alpha+\beta}
			\begin{cases}
				e^{\alpha(s-t)}-e^{\alpha s+\beta t},
				& s\le t,\\
				e^{\beta(t-s)}-e^{\alpha s+\beta t},
				& t\le s,
			\end{cases}\\[2mm]
			G_2(s,t)
			&=\frac1{\alpha+\beta}
			\begin{cases}
				e^{\beta(s-t)}-e^{\beta s+\alpha t},
				& s\le t,\\
				e^{\alpha(t-s)}-e^{\beta s+\alpha t},
				& t\le s.
			\end{cases}
		\end{aligned}
	\end{equation*}
	Here $s,t\le0$ and $G_i(s,t)>0$ for $s,t<0$, $i=1,2$.
	A direct calculation shows
	\begin{equation}\label{eq:kernel-integra}
		\int_{-\infty}^0G_1(s,t)\,dt=\frac{1-e^{\alpha s}}{\alpha\beta}
		\le\frac1c,\quad
		\int_{-\infty}^0G_2(s,t)\,dt
		=\frac{1-e^{\beta s}}{\alpha\beta}
		\le\frac1c.
	\end{equation}
	Since $0<U,V\le\delta$ on $(-\infty,0]$, the Green representation yields
	\begin{align}
		U(s)
		&=U(0)e^{\alpha s}
		+\int_{-\infty}^0G_1(s,t)V^p(t)\,dt,\label{halfline1}\\
		V(s)
		&=V(0)e^{\beta s}
		+\int_{-\infty}^0G_2(s,t)U^q(t)\,dt.\label{halfline2}
	\end{align}

We first show that $U(s)=O(e^{\sigma s})$
as $s\to-\infty$.
For $\rho>0$, a direct calculation shows
\begin{align}
	\int_{-\infty}^0G_1(s,t)e^{\rho t}\,dt
	&=\frac{e^{\rho s}-e^{\alpha s}}
	{(\alpha-\rho)(\beta+\rho)},
	\qquad \rho\ne\alpha,
	\label{powerGreen1}\\
	\int_{-\infty}^0G_2(s,t)e^{\rho t}\,dt
	&=\frac{e^{\rho s}-e^{\beta s}}
	{(\beta-\rho)(\alpha+\rho)},
	\qquad \rho\ne\beta.
	\label{powerGreen2}
\end{align}
Set $V_0(s)=V(0)e^{\beta s}$. By \eqref{halfline2}, we have $0<V_0(s)\le V(s)\le\delta$. 
Since $pq>1$ and $q>1$, by the mean value theorem, we obtain
\[
\begin{aligned}
	0\le V^p(s)-V_0^p(s)\le pq\delta^{p-1/q}
	\bigl(V^{1/q}(s)-V_0^{1/q}(s)\bigr)\le pq\delta^{p-1/q}
	\left(
	\int_{-\infty}^0G_2(s,t)U^q(t)\,dt
	\right)^{1/q}.
\end{aligned}
\]
Since $p\beta=\alpha+2>\alpha$, \eqref{powerGreen1} yields
\[
U(0)e^{\alpha s}
+\int_{-\infty}^0G_1(s,t)V_0^p(t)\,dt
\le C_0e^{\alpha s}
\]
for some $C_0>0$. Consequently, by the
nonnegativity of $G_1$, we obtain
\begin{equation}\label{eq:initial-green-estimate}
	\begin{aligned}
		U(s)&=U(0)e^{\alpha s}
		+\int_{-\infty}^0G_1(s,t)V_0^p(t)\,dt+\int_{-\infty}^0G_1(s,t)(V^p(t)-V_0^p(t))\,dt\\	
		&\le C_0e^{\alpha s}
		+pq\delta^{p-1/q}
		\int_{-\infty}^0G_1(s,t)
		\left(
		\int_{-\infty}^0G_2(t,\tau)U^q(\tau)\,d\tau
		\right)^{1/q}dt.
	\end{aligned}
\end{equation}

For $\varepsilon>0$, set
\[
M_\varepsilon
=\sup_{s\le0}
\frac{U(s)}{e^{\sigma s}+\varepsilon}.
\]
Since $U\le\delta$, we have $M_\varepsilon<\infty$.
By Minkowski's inequality, \eqref{eq:kernel-integra}
and \eqref{powerGreen2}, we obtain
\[
\begin{aligned}
	\left(
	\int_{-\infty}^0G_2(t,\tau)U^q(\tau)\,d\tau
	\right)^{1/q}&\le M_\varepsilon
	\left[
	\left(
	\int_{-\infty}^0G_2(t,\tau)
	e^{q\sigma\tau}\,d\tau
	\right)^{1/q}
	+\varepsilon
	\left(
	\int_{-\infty}^0G_2(t,\tau)\,d\tau
	\right)^{1/q}
	\right]\\
	&\le M_\varepsilon
	\left[
	\frac{e^{\sigma t}}
	{[(\beta-q\sigma)(\alpha+q\sigma)]^{1/q}}
	+c^{-1/q}\varepsilon
	\right].
\end{aligned}
\]
Substituting these estimates into
\eqref{eq:initial-green-estimate} and using \eqref{eq:kernel-integra} and
\eqref{powerGreen1}, we obtain
\[
\begin{aligned}
	U(s)
	&\le C_0e^{\alpha s}
	+pq\delta^{p-1/q}M_\varepsilon
	\left(
	C_\sigma e^{\sigma s}
	+c^{-1-1/q}\varepsilon
	\right)\\
	&\le
	C_0e^{\sigma s}
	+2^{-1}M_\varepsilon(e^{\sigma s}+\varepsilon),
\end{aligned}
\]
which implies
$
M_\varepsilon\le2C_0$. Letting
$\varepsilon\to0$, we obtain $
U(s)=O(e^{\sigma s})$ as $s\to-\infty$.

Since $q\sigma<\beta$, \eqref{halfline2} and
\eqref{powerGreen2} imply
$V(s)=O(e^{q\sigma s})$ as $s\to-\infty$.
Then \eqref{halfline1} yields $U(s)=O(e^{\alpha s})$ since $pq\sigma>\alpha$.
Finally, using $q\alpha=\beta+2>\beta$ in
\eqref{halfline2}, we obtain
$V(s)=O(e^{\beta s})$.
This proves \eqref{sharpdecay}.

By \eqref{EF} and \eqref{sharpdecay}, the corresponding radial pair $(u,v)$ is bounded near the origin.
It follows from Lemma~\ref{Newton}\textup{(iii)} that $(u,v)$
extends across the origin as a positive classical solution
of \eqref{system}.
\end{proof}

\begin{proof}[Proof of Theorem~\ref{classification}]
	By Theorem~\ref{symmetry}, $(u,v)$ is radially
	symmetric. Under the Emden--Fowler transformation
	\eqref{EF}, the corresponding pair $(U,V)$
	is a positive global solution of \eqref{ODE}.

	If $(U,V)$ is neither constant nor periodic,
	Proposition~\ref{remaining} and
	Lemma~\ref{removable} imply that $(u,v)$ extends
	to a positive $C^2$ solution in $\R^n$,
	which contradicts the assumption that the singularity at the origin is nonremovable. Thus $(U,V)$ is either constant or periodic.
	
	If $(U,V)$ is constant, then $(U,V)\equiv(A,B)$, which gives $(u(r),v(r))=(Ar^{-\alpha},Br^{-\beta})$.
	
	Otherwise, $(U,V)$ is nonconstant and periodic. 
	The remaining assertions for the periodic case follow directly
	from Proposition~\ref{periodicfamily} and \eqref{EF}, which completes the proof. 
\end{proof}

\subsection{The ground-state orbit}

We identify the ground-state orbit as the unique positive global
orbit of \eqref{ODE} that is neither an equilibrium nor periodic.

\begin{corollary}\label{homoclinicorbit}
	Assume \eqref{critical} holds. Let $(\cU,\cV)$ be the
	positive radial ground state of \eqref{system} in $\R^n$ centered at the origin and
	normalized by $\cU(0)=1$, and define
	\begin{equation}\label{homoclinic}
		G_1(t)=\ee^{\alpha t}\cU(\ee^t),\qquad
		G_2(t)=\ee^{\beta t}\cV(\ee^t).
	\end{equation}
	Then $(G_1,G_2)$ is a positive homoclinic
	solution of \eqref{ODE} with zero energy.
	Moreover, every positive global solution $(U,V)$
	of \eqref{ODE} that is neither constant nor
	periodic satisfies
	\[
	(U(t),V(t))=(G_1(t+\tau),G_2(t+\tau))
	\]
	for some $\tau\in\R$.
\end{corollary}
\begin{proof}
	Since $(\cU,\cV)$ is a positive radial ground state of\eqref{system} in $\R^n$, $(G_1,G_2)$ is a positive
	global solution of \eqref{ODE} satisfying
	$
	(G_1(t),G_2(t))\to(0,0)
	$
	as $t\to-\infty$.
	Thus it is neither constant nor periodic.
	By Proposition~\ref{remaining}, $(G_1,G_2)$
	is homoclinic to zero and has zero energy.
	
	Now let $(U,V)$ be a positive global solution of
	\eqref{ODE} that is neither constant nor periodic.
	By Proposition~\ref{remaining},
	$(U,V)\to(0,0)$ as $t\to-\infty$.
It follows from	Lemma~\ref{removable} that the
	corresponding radial functions $u$ and $v$
	extend to positive $C^2$ solutions of
	\eqref{system} in $\R^n$.
	By \cite[Theorem~1.1]{GHPY}, there exist
	$a\in\R^n$ and $\lambda>0$ such that
	\[
	u(x)=\lambda^\alpha\cU(\lambda(x-a)),
	\qquad
	v(x)=\lambda^\beta\cV(\lambda(x-a)).
	\]
Since $u$ is radial about the origin and
$\cU$ is strictly radially decreasing,
we have $a=0$.
	Thus, by \eqref{EF} and \eqref{homoclinic}, we have
	$
	(U(t),V(t))=(G_1(t+\log\lambda),G_2(t+\log\lambda)).
	$
\end{proof}

\subsection{Annular critical energy}
We show that removable and nonremovable radial solutions exhibit
different behavior in their annular critical energy. For a radial solution $(u,v)$ of \eqref{system},
define the critical energy on annuli by
\[
\mathcal A(R,L)
=\int_{R<|x|<LR}(u^{q+1}+v^{p+1})\dd x,
\qquad R>0,\quad L>1.
\]
By \eqref{critical} and \eqref{EF}, we obtain
\begin{equation}\label{annularidentity}
	\mathcal A(R,L)
	=\om\int_{\log R}^{\log R+\log L}
	(U^{q+1}+V^{p+1})\dd t.
\end{equation}

\begin{corollary}\label{annularcorollary}
	Assume \eqref{critical} holds. Let $(u,v)$ be a positive radial solution of
	\eqref{system}. Then the following hold:
	\begin{enumerate}
		\item[(i)]
		If the singularity at the origin is nonremovable,
		then there exist constants $0<m\le M<\infty$
		such that
		\begin{equation}\label{eq:annular}
		m\log L\le \mathcal A(R,L)\le M\log L
		\qquad\forall R>0,\ L>1.
		\end{equation}
		In particular, if the Emden--Fowler transform
		$(U,V)$ of $(u,v)$ is nonconstant and periodic
		with least period $T$, then
		\[
		\mathcal A(e^T R,L)=\mathcal A(R,L)
		\qquad\forall R>0,\ L>1.
		\]
		Moreover,
		\begin{equation}\label{eq:annular-asy}
		\frac{\mathcal A(R,L)}{\log L}
		\rightarrow
		\frac{\om}{T}\int_0^T
		(U^{q+1}+V^{p+1})\dd t
		\qquad\text{as }L\to\infty,
		\end{equation}
		uniformly with respect to $R>0$.

		\item[(ii)]
		If $(u,v)$ extends to a positive ground state
		solution in $\R^n$, then
		\[
		\lim_{R\to0}\mathcal A(R,2)
		=\lim_{R\to\infty}\mathcal A(R,2)=0.
		\]
	\end{enumerate}
\end{corollary}

\begin{proof}
	If the singularity at the origin is nonremovable,
	Theorem~\ref{classification} implies that $(U,V)$
	is either constant or periodic. In either case,
	$U^{q+1}+V^{p+1}$ is bounded above and away
	from zero, and hence \eqref{annularidentity}
	gives \eqref{eq:annular}.
	
For a nonconstant periodic solution $(U,V)$
with least period $T$, \eqref{annularidentity}
and periodicity yield
$\mathcal A(e^T R,L)=\mathcal A(R,L)$.

Setting
\[
F(t)=U^{q+1}(t)+V^{p+1}(t),
\qquad
\overline F=\frac1T\int_0^T F(t)\dd t,
\]
we have
\[
\int_0^T(F(t)-\overline F)\dd t=0,
\]
Since $F$ is $T$-periodic, the function
\[
G(t)=\int_0^t(F(s)-\overline F)\dd s
\]
is $T$-periodic and hence bounded. It follows that
\[
\int_a^{a+s}F(t)\dd t
=s\overline F+G(a+s)-G(a)
=s\overline F+O(1),
\]
where $O(1)$ is uniform in $a\in\R$ and $s>0$. Consequently, we obtain
\[
\frac{\mathcal A(R,L)}{\log L}
=\frac{\om}{T}\int_0^T
(U^{q+1}+V^{p+1})\dd t
+O\left(\frac1{\log L}\right),
\]
which implies \eqref{eq:annular-asy}.

Finally, $(ii)$ follows from
$(u,v)\in L^{q+1}(\R^n)\times L^{p+1}(\R^n)$.
\end{proof}

\begingroup\small
\providecommand{\bysame}{\leavevmode\hbox to3em{\hrulefill}\thinspace}
\providecommand{\MR}{\relax\ifhmode\unskip\space\fi MR }
\providecommand{\MRhref}[2]{%
  \href{http://www.ams.org/mathscinet-getitem?mr=#1}{#2}
}
\section*{Acknowledgements}

The first author and third author are supported by the National Natural Science Foundation of
China (NSFC, No.~12671132) and the Science and Technology Commission of
Shanghai Municipality (No.~22DZ2229014). The AI tool ChatGPT was used in the preparation of this manuscript.
The authors assume full responsibility for the mathematical validity,
accuracy, and integrity of the proofs presented in the manuscript.

\providecommand{\href}[2]{#2}

\endgroup
\end{document}